\documentclass{amsart}
\usepackage{etex} 

\usepackage{amssymb}
\usepackage{amsmath}

\usepackage{url}
\usepackage{enumerate}
\usepackage{graphicx}

\usepackage{davide_sty}

\newcommand{\erase}[1]{}

\newtheorem{theorem}{Theorem}[subsection]
\newtheorem{lemma}[theorem]{Lemma}
\newtheorem{proposition}[theorem]{Proposition}
\newtheorem{corollary}[theorem]{Corollary}

\def\noproof{\hfill \rule{3pt}{6pt}}

\newtheorem{_algorithm}[theorem]{Algorithm}
\newenvironment{algorithm}{\begin{_algorithm}\rm}{\noproof\end{_algorithm}}

\newtheorem{_criterion}[theorem]{Criterion}
\newenvironment{criterion}{\begin{_criterion}\rm}{\noproof\end{_criterion}}

\newtheorem{_procedure}[theorem]{Procedure}

\newtheorem{_definition}[theorem]{Definition}
\newenvironment{definition}{\begin{_definition}\rm}{\end{_definition}}

\newtheorem{_remark}[theorem]{\it Remark}
\newenvironment{remark}{\begin{_remark}\rm}{\end{_remark}}

\newtheorem{_example}[theorem]{Example}

\newtheorem{_assumption}[theorem]{Assumption}

\newtheorem{_construction}[theorem]{Construction}

\newtheorem{_claim}[theorem]{Claim}

\newtheorem{_conjecture}[theorem]{Conjecture}

\numberwithin{equation}{section}
\numberwithin{table}{section}
\numberwithin{figure}{section}
\renewcommand{\qed}{\hfill {$\Box$}}

\newcommand{\C}{\mathord{\mathbb C}}

\newcommand{\Q}{\mathord{\mathbb  Q}}
\newcommand{\R}{\mathord{\mathbb R}}

\newcommand{\Z}{\mathord{\mathbb Z}}

\newcommand{\AAA}{\mathord{\mathcal A}}

\newcommand{\DDD}{\mathord{\mathcal D}}

\newcommand{\PPP}{\mathord{\mathcal P}}

\newcommand{\RRR}{\mathord{\mathcal R}}

\newcommand{\VVV}{\mathord{\mathcal V}}

\newcommand{\mapdownsurj}{
\hbox{$\bigm\downarrow$}
\llap{\hbox{\raise 2pt\hbox{$\bigm\downarrow$}}}%
\vstrechmapdown
}

\newcommand{\mapupsurj}{
\hbox{$\bigm\uparrow$}
\llap{\hbox{\raise 2pt\hbox{$\bigm\uparrow$}}}%
\vstrechmapup
}

\newcommand{\inj}{\hookrightarrow}
\newcommand{\surj}{\mathbin{\to \hskip -7pt \to}}

\newcommand{\shortset}[2]{\{ {#1} \,|\, {#2}   \}}

\newcommand{\gen}[1]{\langle {#1}  \rangle}

\newcommand{\tensor}{\otimes}

\newcommand{\sprime}{\sp\prime}

\newcommand{\spar}[1]{\sp{(#1)}}

\newcommand{\sperp}{\sp{\perp}}

\newcommand{\dual}{\sp{\vee}}

\newcommand{\inv}{\sp{-1}}

\newcommand{\Hom}{\mathord{\mathrm{Hom}}}

\newcommand{\OG}{\mathord{\mathrm{O}}}

\newcommand{\id}{\mathord{\mathrm{id}}}

\newcommand{\Aut}{\operatorname{\mathrm{Aut}}\nolimits}

\newcommand{\pr}{\mathord{\mathrm{pr}}}

\newcommand{\rank}{\operatorname{\mathrm{rank}}\nolimits}

\newcommand{\closure}[1]{\overline{#1}}

\newcommand{\weyl}{\mathord{\bf w}}
\newcommand{\iotaX}{\iota_X}
\newcommand{\iotaS}{\iota_S}
\newcommand{\iotaY}{\iota_Y}
\newcommand{\enr}{\varepsilon}
\newcommand{\pistar}{\pi^{*}}

\newcommand{\aut}{\mathord{\rm aut}}

\newcommand{\intf}[1]{\langle #1\rangle}

\newcommand{\YIV}[1]{Y_{\hbox to 2.6mm {\scriptsize\rm IV}, \hskip 1pt #1}}

\newcommand{\OGp}{\OG^{+}}

\newcommand{\Xinduced}{\iota_X^* \RRR_{26}\sperp}
\newcommand{\Yinduced}{\iota_Y^* \RRR_{26}\sperp}
\newcommand{\induced}{\iota^* \RRR_{26}\sperp}

\newcommand{\BStype}[1]{\text{\tt #1}}
\newcommand{\NKtype}[1]{\text{\normalfont\rom{#1}}}

\newcommand{\ADE}{\mathord{\rm ADE}}

\newcommand{\Km}{\mathord{\rm Km}}

\makeatletter
\newcommand*{\rom}[1]{\expandafter\@slowromancap\romannumeral #1@}
\makeatother

\begin{document}
\title[Enriques involutions on singular K3 surfaces]
{Enriques involutions on singular K3 surfaces \\ of small discriminants}
\author{Ichiro Shimada}
\address[I. Shimada]{Department of Mathematics, 
Graduate School of Science, 
Hiroshima University,
1-3-1 Kagamiyama, 
Higashi-Hiroshima, 
739-8526 (Japan)}
\email{ichiro-shimada@hiroshima-u.ac.jp}
\author{Davide Cesare Veniani }
\address[D. C. Veniani]{Institut für Mathematik, 
FB 08 - Physik, Mathematik und Informatik,
Johannes Gutenberg-Universit\"at,
Staudingerweg 9, 4. OG,
55128 Mainz (Germany)}
\email{veniani@uni-mainz.de}

\thanks{The first author was supported by JSPS KAKENHI Grant Number 15H05738, ~16H03926,  and~16K13749. The second author was supported by SFB/TRR45.}

\begin{abstract}
We classify Enriques involutions on a K3 surface, 
up to  conjugation in the automorphism group, in terms of lattice theory. 
We enumerate such involutions on singular K3 surfaces with transcendental lattice of discriminant smaller than or equal to $36$. 
For 11 of these K3 surfaces, we apply Borcherds method to compute the automorphism group of the Enriques surfaces covered by them.
In particular,
we investigate the structure of the two most algebraic Enriques surfaces.
\end{abstract}

\maketitle

\section{Introduction}
Let $X$ be a complex K3 surface. 
We denote by $S_X=H^2(X, \Z)\cap H^{1,1}(X)$ the lattice 
of numerical equivalence classes of divisors on $X$,
and by $T_X$ the orthogonal complement of $S_X$ in $H^2(X, \Z)$,
which we call the \emph{transcendental lattice} of $X$.
Suppose that $X$ is \emph{singular},
that is, the Picard number $\rank S_X$ attains the possible maximum $h^{1,1}(X)=20$.
The \emph{discriminant} of a singular K3 surface $X$
is the determinant of a Gram matrix of $T_X$.
Since $T_X$ is an even positive definite lattice of rank $2$,
the discriminant $d$ of $X$ is a positive integer satisfying $d\equiv 0\; \textrm{or}\; 3 \bmod 4$.
Note that $T_X$ is naturally oriented by the Hodge structure.
By the classical work of Shioda--Inose~\cite{ShiodaInose},
we know that the isomorphism class of 
the oriented lattice $T_X$ determines $X$ up to $\IC$-isomorphism.
\par
An involution $\tilde \enr\colon X \rightarrow X$ 
of a K3 surface $X$ is called an \emph{Enriques involution} 
if $\tilde \enr$ acts freely on $X$. 
Sertöz~\cite{Sertoz2005} gave a simple criterion to determine
whether a singular K3 surface has 
an Enriques involution or not (see Theorem~\ref{thm:sertoez} and also Lee~\cite{Lee2012}).
On the other hand, Ohashi~\cite{ohashi}
showed that each complex K3 surface $X$ 
(not necessarily singular)
has only finitely many Enriques involutions 
up to conjugation in the automorphism group $\Aut(X)$ of $X$,
and that there exists no universal bound for the number of
conjugacy classes of Enriques involutions.
\par
In this paper,
we classify, up to conjugation in $\Aut(X)$, 
all Enriques involutions~$\tilde \enr$ on a singular K3 surface $X$
whose discriminant $d$ satisfies $d\le 36$.
The classification is given in Table~\ref{tab:enriques_quotients}.
We then investigate the automorphism group $\Aut(Y)$
of some of the Enriques surfaces 
$Y=X/\gen{\tilde \enr}$ covered by singular K3 surfaces.
The result is given in Theorem~\ref{thm:ichiro-main} and Table~\ref{table:enrtable}.
\par
As a corollary, we obtain the following.
For $d= 3$, $4$ or $7$,
there exists exactly one singular K3 surface $X_d$ of discriminant $d$
up to $\C$-isomorphism. 
The K3 surfaces $X_3,X_4$, also known as ``the two most algebraic K3 surfaces'', were studied by Vinberg~\cite{Vinberg1983}. 
By Sertöz~\cite{Sertoz2005}, 
neither $X_3$ nor $X_4$ admits any Enriques involution, but $X_7$ does. 
Hence,
following Vinberg, 
we call the Enriques surfaces covered by~$X_7$ 
the \emph{most algebraic Enriques surfaces}.
\begin{theorem}\label{thm:X7}
The singular K3 surface $X_7$ of discriminant $7$
has exactly two Enriques involutions $\tilde \enr_\NKtype{1}$ and $\tilde \enr_\NKtype{2}$ 
up to conjugation in $\Aut(X_7)$.
Let $Y_\NKtype{1}$ and $Y_\NKtype{2}$ be
the quotient Enriques surfaces corresponding to $\tilde \enr_\NKtype{1}$ and $\tilde \enr_\NKtype{2}$,
respectively. 
Then $\Aut(Y_\NKtype{1})$ is finite of order $8$,
and $\Aut(Y_\NKtype{2})$ is finite of order $24$.
\end{theorem}
Nikulin~\cite{Nikulin1984} and Kondo~\cite{Kondo1986} classified all complex Enriques surfaces
whose automorphism group is finite.
It turns out that these Enriques surfaces are
divided into $7$ classes $\NKtype{1}$, $\NKtype{2}$, 
\dots, $\NKtype{7}$,
which we call \emph{Nikulin-Kondo type}.
See~Kondo~\cite{Kondo1986}
for the properties of these Enriques surfaces.
\begin{corollary}\label{cor:X7}
The most algebraic Enriques surfaces
have finite automorphism groups
and their Nikulin-Kondo types are $\NKtype{1}$ and $\NKtype{2}$.
\end{corollary}
In Section~\ref{sec:themostalgebraic} of this paper,
we give explicit models of the most algebraic Enriques surfaces $Y_\NKtype{1}$ and $Y_\NKtype{2}$ as Enriques sextic surfaces.
\begin{remark}
The Néron--Severi lattice and the automorphism group of $X_7$ were 
determined by Ujikawa~\cite{Ujikawa2013}. 
Elliptic fibrations on $X_7$ were studied by Harrache--Lecacheux~\cite{harrache-lecacheux} and 
Lecacheux~\cite{lecacheux}.
\end{remark}
\begin{remark}
Mukai~\cite{Mukai} also realized that $X_7$
has Enriques involutions 
that produce Enriques surfaces of Nikulin-Kondo type $\NKtype{1}$ and $\NKtype{2}$.
\end{remark}
Ohashi~\cite{ohashi} gave a lattice theoretic method to enumerate 
Enriques involutions on certain K3 surfaces.
He then classified in~\cite{ohashi2} all Enriques involutions on 
the Kummer surface $\Km(\Jac(C))$ associated with the jacobian variety
of a generic curve~$C$ of genus~$2$.
We refine and generalize Ohashi's method. Our main result, namely Theorem~\ref{thm:davide-main}, applies to any K3 surface, and we use it in the case of singular K3 surfaces to compile Table~\ref{tab:enriques_quotients}.

For some K3 surfaces $X$, 
the group $\Aut(X)$ can be calculated by 
Borcherds method~(\cite{Bor1},~\cite{Bor2});
for instance, Kondo~\cite{KondoKmJacC} implemented it in order to compute $\Aut(\Km(\Jac(C)))$. 
We apply Borcherds method in order to calculate the automorphism group of some of 
singular K3 surfaces $X$,
and to write the action of $\Aut(X)$ on the nef chamber of $X$ explicitly.
Building on this data, we enumerate all Enriques involutions up to conjugation,
and, using also a result of the preprint~\cite{BS} (see Section~\ref{subsec:LtttoLts}), 
we calculate the automorphism group of the Enriques surfaces covered by these K3 surfaces.

Note that the enumeration of Enriques involutions by Ohashi's method
and by Borcherds method are carried out independently.
The results are, of course, consistent.
We hope that these methods will be applied to many other K3 surfaces 
(with smaller Picard number)
and Enriques surfaces covered by them,
and that in these works,
our general results on a K3 surface admitting an Enriques involution
(Lemma~\ref{lem:surjectivity} and Proposition~\ref{prop:surjectivity})
will be useful.

Recently, many studies on the automorphism groups $\Aut(Y)$
of Enriques surfaces~$Y$ have 
appeared~(\cite{AllcockDolgachev2018},~\cite{MukaiOhashi2015},~\cite{Shimada2017}).
Our result gives a description of $\Aut(Y)$
in terms of its action on the lattice $S_Y$
of numerical equivalence classes of divisors on $Y$.
We expect that this description 
is helpful in the search for a more geometric description of $\Aut(Y)$,
that is, for writing elements of $\Aut(Y)$
as birational self-maps on some projective model of $Y$.
\par
This paper is organized as follows.
In Section~\ref{sec:Preliminaries},
we recall basic facts about lattices, K3 surfaces and Enriques surfaces,
and fix notions and notation.
In Section~\ref{sec:number-of-Enriques-quotients},
we classify all Enriques involutions on singular K3 surfaces 
with discriminant $\le 36$ by
a generalization of Ohashi's method.
In Section~\ref{sec:AutSingK3},
we recall Borcherds method,
and apply it to $11$ singular K3 surfaces
whose transcendental lattices are listed in Table~\ref{table:Table1}.
Recently, many geometric studies of singular K3 surfaces of small discriminant
have appeared (see, for example, ~\cite{women2015},~\cite{harrache-lecacheux},~\cite{lecacheux},~\cite{Utsumi2016}).
We summarize the computational data 
for these $11$ singular K3 surfaces in Table~\ref{tab:D0data}.
In Section~\ref{sec:EnriquesInvols},
we explain an algorithm to calculate 
Enriques involutions and the automorphism groups of the Enriques surfaces
from the data obtained by Borcherds method,
and apply this method to the $11$ singular K3 surfaces.
In Section~\ref{sec:themostalgebraic},
we study the most algebraic Enriques surfaces $Y_\NKtype{1}$ and $Y_\NKtype{2}$.
\par
For the computation,
the first author used {\tt GAP}~\cite{GAP}.
On the web page~\cite{thecompdata},
the computational data concerned with Borcherds method is given explicitly. The second author used {\tt GAP} and {\tt sage} on SageMath~\cite{sagemath}. The computational data concerned with Ohashi's method is available upon request.

\section{Preliminaries} \label{sec:Preliminaries}

\subsection{Lattices}
A \emph{lattice} is a free $\IZ$-module $L$ of finite rank with a $\IZ$-valued non-degenerate symmetric form $\langle\, \,,\,\rangle$. 
The \emph{determinant} $\det L$ of $L$ is the determinant of any Gram matrix of $L$. 
A lattice $L$ is \emph{unimodular} if $\det L=\pm 1$.
A lattice with the same underlying $\IZ$-module as $L$ and symmetric form $n\cdot\langle\, \,,\,\rangle$ is denoted by $L(n)$.
The group of isometries of $L$ is denoted $\OG(L)$. 
We let $\OG(L)$ act on $L$ \emph{from the right}. 
A vector $v$ of a lattice $L$ is called an \emph{$n$-vector} if $\langle v,v\rangle = n$. 
We denote by $\cR_L$ the set of $(-2)$-vectors of a lattice $L$.

A lattice $L$ is \emph{even} if $\langle v,v \rangle \in 2\IZ$ for all $v \in L$; otherwise, it is \emph{odd}. 
The \emph{signature} of a lattice $L$ is the signature of $L \otimes \IR$. 
Analogously, we say that $L$ is \emph{positive definite}, \emph{negative definite} or 
\emph{indefinite} if $L \otimes \IR$ is.
A lattice $L$ of rank $n > 1$ is \emph{hyperbolic} if the signature  is $(1,n-1)$. 
A \emph{positive cone} of a hyperbolic lattice $L$
is one of the two connected components
of $\shortset{v\in L\tensor\R}{\intf{v, v}>0}$.
For  a hyperbolic lattice $L$ and a positive cone $\cP_L$  of $L$,
we denote by $\OG(L,\cP_L)$ the group of isometries of $L$ that preserves $\cP_L$.

The standard positive definite lattices 
associated to Dynkin graphs will be 
denoted $A_n$ ($n \geq 1$), $D_n$ ($n \geq 4)$, $E_6$, $E_7$, $E_8$.

\subsection{Surfaces}
Let $Z$ be a K3 surface or an Enriques surface.
We denote by 
$S_Z$ the lattice of numerical equivalence classes of divisors on $Z$,
and call it the \emph{N\'eron--Severi lattice} of $Z$.
Then $S_Z$ is an even  hyperbolic lattice, provided that $\rank S_Z>1$.
Let $\PPP_Z$ denote the positive cone of $S_Z$   that contains an ample class, 
and let $\RRR_Z$ be the set of $(-2)$-vectors of $S_Z$.
For simplicity, we denote by $\aut(Z)$ the the image of the natural representation
\begin{equation} \label{eq:nat_rho}
 \rho_Z\colon \Aut(Z) \rightarrow \OG(S_Z, \PPP_Z).
\end{equation}
We put
\[
N_Z:=\shortset{x\in \PPP_Z}{\intf{x, [\Gamma]}\ge 0\;\;\textrm{for all curves $\Gamma$ on $Z$}},
\]
and call it the \emph{nef chamber} of $Z$.
It is obvious that the action of $\aut(Z)$ on $\PPP_Z$ preserves $N_Z$.

\subsection{Finite bilinear and quadratic forms} \label{subsec:finite-forms}

A \emph{finite quadratic form} is a finite abelian group $G$ together with a function $q\colon G \rightarrow \IQ/2\IZ$ which satisfies
\[
 q(n\alpha) = n^2 q(\alpha) \text{ for every $\alpha \in G$ and $n \in \IZ$}
\]
such that the function $b(q)\colon G\times G \rightarrow \IQ/\IZ$ defined by
\[
 (\alpha, \beta) \mapsto \frac{q(\alpha + \beta) - q(\alpha) - q(\beta)}{2}
\]
is a finite symmetric bilinear form.
For the sake of simplicity, we will denote by $q$ also the underlying finite abelian group $G$. 
The \emph{length}, i.e. the minimal number of generators, of $G$ (resp. of the $p$-torsion part of $G$) is denoted by $\ell(G)$ (resp. $\ell_p(G)$). 
A subgroup $\Gamma \subset G$ is called \emph{isotropic} if $q|\Gamma = 0$, 
where $q|\Gamma$ denotes the restriction of $q$ to $\Gamma$. 
Given an isotropic subgroup $\Gamma$, 
the quadratic form $q$ descends to the quotient group $\Gamma^\perp/\Gamma$,
where 
\[
\Gamma^\perp:=\shortset{\alpha \in G}{b(q)(\alpha, \gamma)=0 \text{ for every $\gamma \in \Gamma$}};
\]
we denote the resulting finite quadratic form by $q|{\Gamma^\perp/\Gamma}$. 

If $L$ is a lattice, then the group $L^\vee/L$, where $L^\vee := \Hom(L,\IZ) \subset L \otimes \IQ$, is a finite abelian group of order $|\det L|$.
The \emph{discriminant bilinear form} of a lattice $L$ is the finite symmetric bilinear form induced by $\langle\, \,,\,\rangle$
\[
 b(L)\colon L^\vee/L \times L^\vee/L \rightarrow \IQ/\IZ.
\]
If $L$ is even, the \emph{discriminant quadratic form} of $L$ is the finite quadratic form induced by $\langle\, \,,\,\rangle$
\[
 q(L)\colon L^\vee/L \rightarrow \IQ/2\IZ.
\]

Let $\OG(q(L))$ denote the automorphism group of the finite quadratic form $q(L)$, which we let act on $q(L)$ \emph{from the right}.
There is a natural homomorphism
\[
 \OG(L) \rightarrow \OG(q(L)), \qquad g \mapsto q(g).
\]

Let $C_n(e)$ be the cyclic group of order $n$ generated by $e$. For $k\geq 1$, we denote by $u_k$ (resp. $v_k$) the finite quadratic form with underlying group $C_{2^k}(e) \times C_{2^k}(f)$ such that
$\langle e,e \rangle = \langle f,f\rangle = 0$ (resp. $\langle e,e \rangle = \langle f,f\rangle = 1$) and $\langle e,f\rangle = \tfrac{1}{2^k}$.
For $a,b \in \IZ$ prime to each other, we denote by $\langle \tfrac{a}{b} \rangle$ the finite quadratic form with underlying group $C_b(e)$ such that $\langle e,e \rangle = \tfrac{a}{b}$.

\subsection{Genera} Given a pair of non-negative integers $(s_+,s_-)$ and a non-degenerate finite quadratic (resp. bilinear) form $h$, the \emph{genus} $\fg(s_+,s_-,h)$ is the set of isometry classes of even (resp. odd) lattices of signature $(s_+,s_-)$ with discriminant quadratic (resp. bilinear) form isomorphic to~$h$. 
If a genus contains only the isometry class of a lattice $L$, we say that \emph{$L$ is unique in its genus}.

In general, enumerating all isometry classes in a given genus is a non-trivial problem. It is computationally easier to find lattices of smaller determinant, so the following elementary lemma can be very useful.

\begin{lemma} \label{lem:divisible-lattice}
Given a lattice $L$ and a prime number $p$, then $\ell_p(L^\vee/L) = \rank L$
if and only if $L = L'(p)$ for some lattice $L'$.
In this case, and if moreover $L$ is even and $p = 2$, then $L'$ is odd if and only if 
$q(L) = \anglefrac{1}{2} \oplus q'$ or $q(L) = \anglefrac{3}{2} \oplus q'$
for some finite quadratic form $q'$.
\noproof
\end{lemma}

\begin{remark}
Suppose $q$ is a finite quadratic form admitting an isotropic subgroup~$\Gamma$. In order to enumerate all isometry classes of even lattices in $\fg(s_+,s_-,q)$, we can take advantage of Proposition 1.4.1 in~\cite{Nikulin1979}: first we enumerate all lattices in $\fg(s_+,s_-,q|\Gamma^\perp/\Gamma)$, then we inspect all sublattices of index $|\Gamma|$.
\end{remark}

Given a finite (bilinear or quadratic) form $h$ and $s \in \IN$, the following algorithm, suggested by Degtyarev~\cite{Degtyarev}, finds all (odd or even) lattices in $\fg(s,0,h)$. If $h$ is quadratic we put $b = b(h)$, otherwise we put $b = h$.

\begin{algorithm} \label{alg:lattices-in-genus} 
 Let $r$ be the smallest possible rank for which there exists an (odd or even) positive definite lattice $M$ of rank $r$ and discriminant bilinear form~$-b$.
 By results of Nikulin~\cite{Nikulin1979}, for each $N \in \fg(s,0,h)$ there exists a primitive embedding $\iota\colon M \hookrightarrow L$ into some positive definite unimodular lattice~$L$ of rank $r + s$ such that $[\iota]^\perp \cong N$.
 Taking advantage of the classification of positive definite unimodular lattices of small rank (see, for instance, Table~16.7 in~\cite{conway-sloane}), we list all such lattices $L$. 
 Using {\tt GAP} and the function {\tt ShortestVectors}, we list all primitive embeddings $\iota\colon M\hookrightarrow L$ for all $M \in \fg'$ and all $L$. 
 Then, we compute the lattices~$[\iota]^\perp$ and select those ones which belong to $\fg(s,0,h)$.
In order to eliminate pairs of isomorphic lattices, one can use the attribute \verb+is_globally_equivalent_to+ of the class \verb+QuadraticForm+ in {\tt sage}. 
\end{algorithm}

The algorithm works provided that $r + s$ is small enough and that we can find a lattice $M$ explicitly. In order to find $M$, we can apply the algorithm recursively to $\fg(r,0,-b)$. If $r = 1$ or $2$, this genus can be enumerated a priori (see, for instance, Chapter~15 in~\cite{conway-sloane}). 

\begin{remark}
Another well-known way to enumerate lattices in a given genus is Kneser's neighboring method~\cite{Kneser57}. This method has been implemented in {\tt sage} by Brandhorst~\cite{Brandhorst}.
\end{remark}

\subsection{Primitive embeddings} \label{sec:primitive-sublattices}

Given an embedding of lattices $\iota\colon M \hookrightarrow S$, 
we denote by $[\iota]$ its image and by $[\iota]^\perp$ the orthogonal complement of $[\iota]$ in $S$. 
An embedding $\iota\colon M \hookrightarrow S$ is called \emph{primitive} if $S/[\iota]$ is a torsion-free group. 
All primitive embeddings are considered up to the action of $\OG(M)$.

\begin{proposition}[Proposition 1.15.1 in \cite{Nikulin1979}] \label{prop:Nikulin1.15}
If $\iota\colon M \hookrightarrow S$ is a primitive embedding of even lattices, then there exist a subgroup $H \subset M^\vee/M$ and an isomorphism of finite quadratic forms $\beta\colon q([\iota])|H \rightarrow q(S)|\beta(H)$ such that
\[
 q([\iota]^\perp) \cong (-q([\iota])) \oplus q(S)|{\Gamma^\perp_\beta/\Gamma_\beta},
\]
where $\Gamma_\beta$ is the push-out of $\beta$ in $(-q([\iota])) \oplus q(S)$.
\noproof
\end{proposition}

Given a primitive embedding $\iota\colon M \hookrightarrow S$, we put
\[
 \OG(S,[\iota]) := \{g \in \OG(S)\,|\,[\iota]^g = [\iota]\},
\]
and we denote by $\OG(q(S),[\iota])$ its image in $\OG(q(S))$ by the natural homomorphism $\OG(S) \rightarrow \OG(q(S))$.

Fix now two even lattices $M$, $N$ and consider the set $I(S,M,N)$ of primitive embeddings $\iota\colon M \hookrightarrow S$ such that $[\iota]^\perp \cong N$. 
The group $\OG(S)$ acts on $I(S,M,N)$ in a natural way. 

Consider 
also the set of pairs $(H,\gamma)$, where $H \subset M^\vee/M$ is a subgroup and $\gamma \colon q(M)|H \rightarrow -q(N)|{\gamma(H)}$ is an isomorphism of finite quadratic forms such that
\begin{equation} \label{eq:Hgamma-pushout}
 q(M) \oplus q(N)|{\Gamma_\gamma^\perp/\Gamma_\gamma} \cong q(S),
\end{equation}
where $\Gamma_\gamma$ is the push-out of $\gamma$ in $q(M) \oplus q(N)$. 
We say that two such pairs $(H,\gamma)$ and $(H',\gamma')$ are \emph{equivalent} 
if there exist $\varphi \in \OG(M)$ and $\psi \in \OG(N)$ such that $H^{q(\varphi)} = H'$ and
\begin{equation} \label{eq:gamma-varphi-phi}
 \gamma' \circ q(\varphi) = q(\psi) \circ \gamma.
\end{equation}

\begin{proposition}[Proposition 1.5.1 in \cite{Nikulin1979}] \label{prop:Nikulin1.5.1}
In the above notation, there is a one-to-one correspondence between the elements of $I(S,M,N)$ modulo the action of $\OG(S)$ and the set of pairs $(H,\gamma)$ modulo equivalence.
\noproof
\end{proposition}

\begin{proposition}[Proposition 1.5.2 in \cite{Nikulin1979}] \label{prop:OMS}
For a fixed pair $(H,\gamma)$ corresponding to the orbit of a primitive embedding $\iota\colon M \hookrightarrow S$, the subgroup $\OG(q(S),[\iota])$ consists of those elements $\xi \in \OG(q(S))$ for which there exist $\varphi \in \OG(M)$ and $\psi \in \OG(N)$ such that $H^{q(\varphi)} = H$, equation \eqref{eq:gamma-varphi-phi} holds, and $\xi$ corresponds under the isomorphism~\eqref{eq:Hgamma-pushout} to the automorphism induced by $\varphi$ and $\psi$ on $\Gamma_\gamma^\perp/\Gamma_\gamma$.
\noproof
\end{proposition}

\subsection{Chambers and their faces}\label{subsec:defchambersfaces}
Let $V$ be a $\Q$-vector space of dimension $n>1$
with a non-degenerate symmetric bilinear form 
$\intf{\;, \;}\colon V\times V\to \Q$
such that $V\tensor \R$ is of signature $(1, n-1)$.
Let $\PPP_V$ be
one of the two connected components
of $\shortset{x\in V\tensor\R}{\intf{x, x}>0}$.
For $v\in V$ with $\intf{v, v}<0$,
we put 
\[
(v)\sperp:=\shortset{x\in \PPP_V}{\intf{x, v}=0},
\]
which is a hyperplane of $\PPP_V$.
For a set $\VVV$ of vectors $v\in V$ with $\intf{v, v}<0$,
we denote by $\VVV\sperp$ the family of hyperplanes
$\shortset{(v)\sperp}{v\in \VVV}$.
\par
Let $\VVV$ be a set of vectors $v\in V$ with $\intf{v, v}<0$ such that 
the family of hyperplanes $\VVV\sperp$ is locally finite.
A \emph{$\VVV\sperp$-chamber} is the closure in $\PPP_V$ of a connected component of 
the complement
\[
\PPP_V\setminus \bigcup_{H\in \VVV\sperp} H.
\]
Let $\closure{\PPP}_V$ be the closure of $\PPP_V$ in $V\tensor \R$,
and $\partial\,\closure{\PPP}_V$ 
the boundary $\closure{\PPP}_V\setminus \PPP_V$ of $\closure{\PPP}_V$.
Let $C$ be a $\VVV\sperp$-chamber,
and $\closure{C}$ the closure of $C$ in $V\tensor \R$.
We say that $C$ is \emph{quasi-finite}
if $\closure{C} \cap \partial\, \closure{\PPP}_V$
is contained in a union of at most countably many real half-lines of $V\tensor \R$.
\par
Let $C$ be a quasi-finite $\VVV\sperp$-chamber.
Suppose that we are given a set $U_C$
of vectors $v\in V$ with $\intf{v, v}<0$ such that
\[
C=\shortset{x\in \PPP_V}{\intf{x, v}\ge 0\;\;\textrm{for all}\;\; v\in U_C}.
\]
A \emph{wall} of $C$ is a closed subset $w$ of $C$ 
for which 
there exists a hyperplane $H\in \VVV\sperp$ with $w=C\cap H$
such that $w$ contains a non-empty open subset of $H$.
Let $w$ be a wall of $C$. 
A vector $v\in V$ with $\intf{v,v}<0$  is said to \emph{define $w$} 
if $w$ is equal to $C\cap (v)\sperp$
and $\intf{x, v}>0$ holds for all interior points $x$ of $C$.
A vector $v_0\in U_C$ defines a wall of $C$ if and only if
there exists a point $y\in \PPP_V$ such that
$\intf{y, v_0}< 0$ and that $\intf{y, v\sprime}> 0$ holds 
for all $v\sprime\in U_C$ with $(v\sprime)\sperp\ne (v_0)\sperp$.
Therefore, if $U_C$ is finite,
we can calculate the set of walls of $C$ by means of linear programming.
\par
A \emph{face} is a closed subset of $C$ that is the intersection of a finite number of walls of $C$.
Let $f$ be a face of $C$.
We denote by $\gen{f}$ the minimal linear subspace of $V$ containing $f$.
The \emph{dimension of $f$} is the dimension of $\gen{f}$.
Suppose that $m:=\dim f$ is $\ge 2$.
Since $f$ contains a non-empty open subset of $\gen{f}$,
the linear space $\gen{f}$ contains a vector $v$ with $\intf{v, v}>0$,
and hence the restriction of $\intf{\;, \;}$ to $\gen{f}$ is of signature $(1, m-1)$.
We denote by 
\[
\iota_{\gen{f}}\colon \gen{f} \inj V\;\; {\rm and}\;\; \pr_{\gen{f}}\colon V\surj \gen{f}
\]
the inclusion and the orthogonal projection, respectively,
and let $\PPP_{\gen{f}}$ be the positive cone of $\gen{f}$ that is mapped into $\PPP_V$ by $\iota_{\gen{f}}$.
We put
\[
\iota_{\gen{f}}^*\VVV\sperp:=\shortset{\iota_{\gen{f}}\inv (H)}{ \textrm{$ H\in \VVV\sperp$ such that $\iota_{\gen{f}}\inv(H)$ is a hyperplane of $\PPP_{\gen{f}}$}},
\]
which is a locally finite family of hyperplanes of $\PPP_{\gen{f}}$.
Note that $\iota_{\gen{f}}^*\VVV\sperp$ is equal to $(\pr_{\gen{f}}^*\VVV)\sperp$,
where
\[
\pr_{\gen{f}}^*\VVV:=\shortset{\pr_{\gen{f}} (v)}{ v \in \VVV\;\;\textrm{such that}\;\; \intf{\pr_{\gen{f}} (v), \pr_{\gen{f}} (v)}<0}.
\]
Then the face $f$ of $C$ is an $\iota_{\gen{f}}^*\VVV\sperp$-chamber in $\PPP_{\gen{f}}$,
and is equal to
\[
\shortset{z\in \PPP_{\gen{f}}}{\textrm{$\intf{z, \pr_{\gen{f}}(v)}\ge 0$ for all $v\in U_C$ with $\intf{\pr_{\gen{f}}(v), \pr_{\gen{f}}(v)}<0$}}.
\]
Therefore, if $U_C$ is finite,
we can calculate the set of walls of the $\iota_{\gen{f}}^*\VVV\sperp$-chamber~$f$,
and hence we can calculate the set of all faces of $C$ by descending induction on the dimension of faces.
\par
Let $w$ be a wall of $C$. Then there exists a unique $\VVV\sperp$-chamber $C\sprime$ such that
$C\cap C\sprime=w$.
This $\VVV\sperp$-chamber $C\sprime$ is said to be \emph{adjacent to $C$ across the wall $w$}.
\subsection{Induced chambers}
Let $L$ be an even hyperbolic lattice.
We apply the above definitions to $L\tensor \Q$.
Let $\PPP_L$ be a positive cone of $L$,
and let $\VVV$ be a set of vectors $v \in L\tensor \Q$ with $\intf{v, v}<0$ such that 
the family $\VVV\sperp$ of hyperplanes of $\PPP_L$ is locally finite.
Suppose that we have a primitive embedding
\[
\iotaS\colon S\inj L
\]
of an even hyperbolic lattice $S$ of rank $m<n$,
and let $\PPP_S$ be the positive cone of $S$
that is mapped into $\PPP_L$ by $\iotaS$.
We use the same letter $\iotaS$ to denote the inclusion $\PPP_S\inj \PPP_L$.
We denote the orthogonal projection by
$\pr_S\colon L\tensor \Q\to S\tensor \Q$,
and put 
\begin{align*}
\iotaS^*\VVV\sperp&:= \shortset{\iotaS\inv (H)}{ \textrm{$ H\in \VVV\sperp$ 
such that $\iotaS\inv(H)$ is a hyperplane of $\PPP_S$}},\\
\pr_S^*\VVV &:= \shortset{\pr_S(v)}{ \text{$v\in \VVV$ with $\intf{\pr_S(v), \pr_S(v)}<0$}}.
\end{align*}
Then $\iota_S^*\VVV\sperp=(\pr_S^*\VVV)\sperp$ is a locally finite family of hyperplanes of $\PPP_S$.
A $\VVV\sperp$-chamber $C\subset \PPP_L$ is said to be \emph{non-degenerate} with respect to $\iotaS$
if the closed subset $\iotaS\inv (C)$ of $\PPP_S$ contains a non-empty open subset of $\PPP_S$.
Suppose that $C$ is non-degenerate with respect to $\iotaS$.
Then $\iotaS\inv (C)$ is an $\iotaS^*\VVV\sperp$-chamber, 
which we call the chamber \emph{induced by $C$}.
If $C$ is quasi-finite, then so is the induced chamber~$\iotaS\inv (C)$.
%
\subsection{Vinberg chambers and Conway chambers}
Let $L$ be as above.
Note that the family $\RRR_L\sperp$ of hyperplanes is locally finite,
where $\RRR_L$ is the set of $(-2)$-vectors.
Each $r\in \RRR_{L}$ defines a reflection
$x\mapsto x+\intf{x, r} r$.
Let $W(L)$ be the subgroup of $\cO(L,\cP_L)$ generated by reflections with respect to $(-2)$-vectors. 
Then each $\RRR_L\sperp$-chamber is a standard fundamental domain of the action of $W(L)$ on $\PPP_L$.

For $n=10$ and $n=26$, 
let $L_{n}$ be an even unimodular hyperbolic lattice of rank~$n$,
which is unique up to isomorphism.
We denote by $\PPP_{n}$ a positive cone of $L_n\tensor \R$, 
and by $\RRR_{n}$ the set of $(-2)$-vectors of $L_n$.

An $\RRR_{10}\sperp$-chamber in $\PPP_{10}$ is called a \emph{Vinberg chamber}.
It is known that a Vinberg chamber is quasi-finite.
\begin{theorem}[Vinberg~\cite{Vinberg1975}]
A Vinberg chamber has exactly $10$ walls.
\noproof
\end{theorem}
An $\RRR_{26}\sperp$-chamber in $\PPP_{26}$ is called a \emph{Conway chamber}.
It is known that a Conway chamber is quasi-finite.
A non-zero primitive vector $\weyl\in L_{26}\cap \, \partial\, \PPP_{26}$ 
is called a \emph{Weyl vector} if 
the negative definite lattice $[\weyl]\sperp/[\weyl]$ is isomorphic to the negative definite Leech lattice,
where $[\weyl]\sperp:=\shortset{v\in L_{26}}{\intf{v, \weyl}=0}$.
\begin{theorem}[Conway~\cite{Vinberg1975}]\label{thm:Conway}
For each Conway chamber $C$,
there exists a unique Weyl vector $\weyl_C$ such that the walls of $C$ are defined by
$(-2)$-vectors $r\in \RRR_{26}$ satisfying $\intf{\weyl, r}=1$.
\noproof
\end{theorem}
\subsection{Primitive embeddings of \texorpdfstring{$L_{10}(2)$}{L10(2)} into \texorpdfstring{$L_{26}$}{L26}} \label{subsec:LtttoLts}
In~\cite{BS},
we classified all primitive embeddings 
of $L_{10}(2)$ into $L_{26}$. 
It turns out that,
up to the action of $\OG(L_{10}(2))=\OG(L_{10})$ and $\OG(L_{26})$,
there exist exactly $17$ primitive embeddings,
which are named as being of type 
\[
\BStype{12A}, \BStype{12B}, \;\BStype{20A}, \dots, \BStype{20A},\dots, \BStype{20F},\;
\BStype{40A}, \dots, \BStype{40E},\; \BStype{96A},
\dots, \BStype{96C},\; \BStype{infty}.
\]
Let $\iota\colon L_{10}(2)\inj L_{26}$
be a primitive embedding.
Identifying 
positive cones of $L_{10}(2)$ with positive cones of $L_{10}$
and replacing $\iota$ with $-\iota$ if necessary, 
we assume that $\iota$ maps $\PPP_{10}$ into $\PPP_{26}$.
Then $\PPP_{10}$ is covered by $\induced$-chambers.
Since Conway chambers are quasi-finite,
every $\induced$-chambers are quasi-finite.
In~\cite{BS},
we have proved the following:
\begin{theorem}\label{thm:BS}
Suppose that $\iota$ is not of type \BStype{infty}.
Let $D$ and $D\sprime$ be $\induced$-chambers. 
Then there exists an isometry $g\in \OGp(L_{10})$
that preserves the set of $\induced$-chambers
and maps $D$ to $D\sprime$.
Each $\induced$-chamber has only a finite number of walls,
and each wall is defined by a $(-2)$-vector.
If $D\cap (r)\sperp$ is a wall of $D$
with $r\in \RRR_{10}$, then
the $\induced$-chamber adjacent to $D$ across the wall $D\cap (r)\sperp$
is the image of the reflection of $D$ into the hyperplane $(r)\sperp$.
\noproof
\end{theorem}
\begin{remark}\label{rem:infty}
If a  primitive embedding $\iota\colon L_{10}(2)\inj L_{26}$ is of type \BStype{infty},
then the $\induced$-chamber has infinitely many walls.
The embedding  $\iota$ is of type \BStype{infty}
if and only if $[\iota]\sperp$ 
contains no $(-2)$-vectors.
\end{remark}
Let $Y$ be an Enriques surface.
Then the N\'eron-Severi lattice $S_Y$ is isomorphic to~$L_{10}$.
It is known that the nef chamber $N_Y$ is bounded by hyperplanes $(r)\sperp$
defined by $(-2)$-vectors $r\in \RRR_Y$.
In~\cite{BS}, we have proved the following:
\begin{theorem}\label{thm:sigmatau}
Let $[\sigma, \tau]$ be one of the pairs
\[
[\BStype{12A}, \NKtype{1}], \;
[\BStype{12B}, \NKtype{2}], \;
[\BStype{20A}, \NKtype{5}], \;
[\BStype{20B}, \NKtype{3}], \;
[\BStype{20C}, \NKtype{7}], \;
[\BStype{20D}, \NKtype{7}], \;
[\BStype{20E}, \NKtype{6}], \;
[\BStype{20F}, \NKtype{4}].
\]
Then every $\induced$-chamber $D$ for a primitive embedding
$\iota\colon L_{10}(2) \inj L_{26}$ of type $\sigma$
is equal to the nef chamber $N_Y$ of 
an Enriques surface $Y$
with finite automorphism group
of Nikulin-Kondo type $\tau$
under an isomorphism $L_{10}\cong S_Y$.
\noproof
\end{theorem}
\subsection{K3 surfaces} 
Let $X$ be a complex projective K3 surface with transcendental lattice~$T_X$. 
Then 
the nef chamber $N_X$ is an $\cR_X^\perp$-chamber, 
and each wall of $N_X$ is defined by the class of a smooth rational curve on $X$. 
We put
\[
 \OG(S_X,N_X) := \{g \in \OG(S_X)\,|\,N_X^g = N_X \}.
\]
Recall that $W_X:=W(S_X)$ is the subgroup of $\OG(S_X, \PPP_X)$ generated by reflections
with respect to $(-2)$-vectors.
The following relations hold (see \cite{ohashi}):
\begin{align} \label{eq:OpS = WxOaS}
& \OG(S_X,\cP_X) = W_X \rtimes  \OG(S_X,N_X), \\
& \label{eq:WinKer} W_X \subset \ker( \OG(S_X) \rightarrow \OG(q(S_X))).
\end{align}

Let $\OG(T_X, \omega_X)$ be the group of isometries of $T_X$ that preserves the $1$-dimensional subspace $H^{2,0}(X) \subset T_X \otimes \IC$, and let $\OG(q(T_X), \omega_X)$ be the image of $\OG(T_X, \omega_X)$ by the natural homomorphism $\OG(T_X) \rightarrow \OG(q(T_X))$. 
The even unimodular overlattice
$H^2(X,\IZ)$ of the orthogonal direct sum $S_X \oplus T_X$ induces an anti-isometry between
the discriminant forms of $S_X$ and of $T_X$ (see \cite{Nikulin1979}), and hence induces an isomorphism $\OG(q(S_X)) \cong \OG(q(T_X))$. 
Let $\OG(q(S_X),\omega_X)$ be the image of $\OG(q(T_X),\omega_X)$ through this isomorphism.
We say that an isometry $g \in \OG(S_X)$ \emph{satisfies the period condition} if $q(g) \in \OG(q(S_X),\omega_X)$. 
Let $\OG(S_X, \omega_X)$ denote
the group of isometries satisfying the period condition.
Recall that $\aut(X)\subset \OG(S_X,\PPP_X)$ is the image of $\Aut(X)$ by~\eqref{eq:nat_rho}.
The Torelli theorem for complex K3 surfaces asserts that
\begin{equation} \label{eq:torelli}
 \aut(X) =  \OG(S_X,N_X) \cap \OG(S_X,\omega_X).
\end{equation}
In particular, if $g\in \OG(S_X, \omega_X)$ 
maps an interior point of $N_X$ to an interior point of $N_X$,
then $g$ belongs to $\aut(X)$.
\begin{remark}\label{rem:Kerrho}
By the Torelli theorem, the kernel of 
$\rho_X\colon \Aut(X)\to \OG(S_X)$
is isomorphic to the kernel of 
the natural homomorphism $\OG(T_X, \omega_X)\to \OG(q(T_X))$.
\end{remark}

\subsection{Singular K3 surfaces} \label{sec:singular-K3-surfaces}
Let $X$ be a singular K3 surface. 
Its transcendental lattice $T_X$ admits a basis with respect to which the Gram matrix is of the form 
\[
    [a,b,c] := \begin{bmatrix} a & b \\ b & c \end{bmatrix},
\] 
with $0 \leq 2b \leq a \leq c$. 
We write $X(T)$ for the K3 surface corresponding to an oriented positive definite even lattice $T$ of rank~$2$.
The lattice $\overline T = [a,-b,c]$ defines a distinct oriented isomorphism class if and only if
$0 < 2b < a < c$.

\begin{remark} \label{rmk:OT-omega-singular-K3}
If $X$ is a singular K3 surface, the subgroup $\OG(T_X,\omega_X)$ can be identified with the subgroup consisting of isometries of $T_X$ of positive determinant. Its image $\OG(q(T_X),\omega_X)$ depends only on the genus of $T_X$. 
\end{remark}


\section{Classification of Enriques involutions up to conjugation} \label{sec:number-of-Enriques-quotients}

Let $X$ be a complex projective K3 surface.
We are interested in classifying the images $\enr$ of Enriques involutions $\tilde\enr$
in $\aut(X)$ 
through the natural representation~\eqref{eq:nat_rho} up to conjugation in $\aut(X)$.
The image $\enr\in \aut(X)$ is also call an Enriques involution.
This is essentially the same problem by the following observation due to Ohashi.

\begin{proposition}[Ohashi~\cite{ohashi}] \label{prop:ohashi-conj}
Let $\tilde \enr_{1},\tilde \enr_2\colon X\rightarrow X$ be two Enriques involutions. 
Then the quotients $Y_i := X/\langle \tilde \enr_i\rangle$, $i = 1,2$, are isomorphic over $\IC$ if and only if $\enr_{1}$, $\enr_2$ are conjugate in $\aut(X)$.
\noproof
\end{proposition}

In this section, after recalling part of Ohashi's work, 
we refine and generalize his main Theorem 2.3 in~\cite{ohashi}.

\subsection{Main result} 
Given an Enriques involution~$\enr \in \aut(X)$, we put
\[
S_X^{\enr=1} := \shortset{v \in S_X}{v^\enr = v}.
\]
We have the following criterion by Keum.
\begin{theorem}[Keum~\cite{Keum1990}] \label{thm:keum}
An involution $\enr\in \aut(X)$ is an Enriques involution 
if and only if the following holds:
the sublattice $S_X^{\enr=1}$
is isomorphic to $L_{10}(2)$ and its orthogonal complement
in $S_X$ contains no $(-2)$-vectors.
\noproof
\end{theorem}

Let $I_X$ be the set of primitive embeddings $\iota\colon L_{10}(2) \hookrightarrow S_X$ such that the orthogonal complement $[\iota]^\perp$ of the image of $\iota$ in $S_X$ contains no $(-2)$-vectors.
The group $\OG(S_X)$ acts on $I_X$ in a natural way.

\begin{proposition}[Proposition 2.2 in \cite{ohashi}] \label{prop:existence-ohashi}
For every $\iota \in I_X$ and $g \in \OG(S_X)$ such that $[\iota]^g$ intersects the interior of $N_X$, there exists a unique $\enr \in \aut(X)$ such that $S_X^{\enr=1} = [\iota]^g$.
\noproof
\end{proposition}

\begin{corollary} \label{cor:easy}
 Let $\enr_1,\enr_2 \in \aut(X)$ be two Enriques involutions.
 Then, there exists $\gamma \in \aut(X)$ such that $\enr_2 = \gamma \circ \enr_1 \circ \gamma^{-1}$ if and only if $(S_X^{\enr_1=1})^\gamma = S_X^{\enr_2=1}$.
 \noproof
\end{corollary}

\begin{proposition}[Step 1 of Theorem 2.3 in \cite{ohashi}] \label{prop:step1}
 For every $\iota \in I_X$ there exists $h \in \OG(S_X)$ such that $[\iota]^h$ intersects the interior of $N_X$.
 \noproof
\end{proposition}

\begin{lemma}[Step 2 of Theorem 2.3 in \cite{ohashi}] \label{lem:w_acts_trivially}
Suppose $[\iota]$ intersects the interior of $N_X$. If there exist an Enriques involution $\enr \in \aut(X)$ and $g \in \OG(S_X)$ such that $S_X^{\enr=1} = [\iota]^g$, 
then there exists $\tilde g \in \OG(S_X,N_X)$ such that $S_X^{\enr=1} = [\iota]^{\tilde g}$.
\noproof
\end{lemma}

\begin{proposition} \label{prop:main-Davide}
Given $\iota \in I_X$, let $\enr_1,\enr_2 \in \aut(X)$ be two Enriques involutions with $S_X^{\enr_1=1} = [\iota]^{g_1}$ and $S_X^{\enr_2=1} = [\iota]^{g_2}$ for some $g_1,g_2 \in \OG(S_X,N_X)$.
Then the Enriques involutions $\enr_1$ and $\enr_2$ are conjugate in $\aut(X)$ if and only if the natural images $q(g_1),q(g_2) \in \OG(q(S_X))$ belong to the same double coset with respect to 
$\OG(q(S_X),[\iota])$ and $\OG(q(S_X),\omega_X)$.
\end{proposition}
\proof
Let $\iota_i := g_i\circ\iota$ for $i = 1,2$. Suppose there exists $\gamma \in \aut(X)$ with $\enr_2 = \gamma \circ \enr_1 \circ \gamma^{-1}$.
Let $\varphi := g_2^{-1}\circ\gamma\circ g_1$, so that $\varphi \in \OG(S_X, [\iota])$. Indeed, by Corollary~\ref{cor:easy},
\[
 [\iota]^\varphi = [\gamma\circ\iota_1]^{g_2^{-1}} = [\iota_2]^{g_2^{-1}} = [\iota].
\]
As $g_1 = \varphi \circ g_2 \circ \gamma^{-1}$ and $\gamma \in \OG(S_X,\omega_X)$, 
the automorphisms $q(g_1)$, $q(g_2)$ of $q(S_X)$ belong to the same double coset.

Conversely, assume that there exist $\varphi \in \OG(S_X,[\iota])$ and $\gamma' \in \OG(S_X,\omega_X)$ such that $q(g_2) = q(\varphi \circ g_1 \circ {\gamma'})$ in $\OG(q(S_X))$. 
Without loss of generality, we can suppose $\varphi \in \OG(S_X,N_X)$. In fact, we can first exchange $\varphi$ with $-\varphi$ if necessary and suppose that $\varphi \in \OG(S_X,\cP_X)$. 
By~\eqref{eq:OpS = WxOaS} and \eqref{eq:WinKer}, we can write $\varphi = w \circ \varphi'$, with $w \in W_X$ and $\varphi'\in \OG(S_X,N_X)$ and exchange $\varphi$ with $\varphi'$ if necessary. 
Define now $\gamma := g_2 \circ \varphi^{-1} \circ g_1^{-1}$. 
Then $\gamma \in \OG(S_X,N_X)$ and $q(\gamma) = q(\gamma')$, so $\gamma \in \OG(S_X,\omega_X)$.
The Torelli Theorem~\eqref{eq:torelli} implies that $\gamma \in \aut(X)$. 
Furthermore, we have
\[
 [\iota_1]^\gamma = ([\iota]^{\varphi^{-1}})^{g_2} = [\iota_2],
\]
so $\enr_1$ and $\enr_2$ are conjugate in $\aut(X)$ by Corollary~\ref{cor:easy}.
\endproof

\begin{lemma} \label{lem:surjectivity}
If a K3 surface $X$ admits at least one Enriques involution, 
then the lattice $S_X$ is unique in its genus and the natural homomorphism $\OG(S_X) \rightarrow \OG(q(S_X))$ is surjective.
\end{lemma}
\proof
Let $\iota\colon L_{10}(2) \hookrightarrow S_X$ be a primitive embedding. 
Then $q(S_X) \cong (q([\iota]) \oplus q([\iota]^\perp))|{\Gamma^\perp/\Gamma}$ for some isotropic subgroup $\Gamma$ of $q([\iota]) \oplus q([\iota]^\perp)$. 
Since $q([\iota]) \cong q(L_{10}(2)) \cong u_{1}^{\oplus 5}$, this implies that 
\[
\ell_p(S_X^\vee/S_X) \leq \rank [\iota]^\perp = \rank S_X - 10
\]
for every odd prime $p$. Moreover, if $\ell_2(S_X^\vee/S_X) = \rank S_X$, then $q(S_X) = q([\iota]) \oplus q'$ for some finite quadratic form $q'$.
Therefore, we can conclude by Theorem 1.14.2 in \cite{Nikulin1979}.
\endproof

Combining Lemma~\ref{lem:surjectivity} and the same argument as in Step 5 of Theorem~2.3 in~\cite{ohashi}, we prove the following proposition.

\begin{proposition} \label{prop:surjectivity}
If a K3 surface $X$ admits at least one Enriques involution, then $\OG(S_X,N_X) \rightarrow \OG(q(S_X))$ is surjective.
\qed
\end{proposition}

Our main result is the following theorem.

\begin{theorem} \label{thm:davide-main}
Let $X$ be a K3 surface and $\iota_{1},\ldots,\iota_r \in I_X$ be a complete set of representatives for the action of $\OG(S_X)$ on $I_X$. Then there exists a bijection between the set of Enriques involutions up to conjugation in $\aut(X)$ and the disjoint union of the sets of double cosets
\[
 \OG(q(S_X),[\iota_i])\backslash\OG(q(S_X))/\OG(q(S_X),\omega_X),\quad i = 1,\ldots,r.
\]
\end{theorem}
\proof
Let $G = \OG(S_X)$, $H_i = \OG(q(S_X),[\iota_i])$ and $K = \OG(q(S_X),\omega_X)$. For each $i = 1,\ldots,r$, fix $h_i \in G$ such that $[\iota_i]^{h_i}$ intersects the interior of $N_X$ (Proposition~\ref{prop:step1}). As exchanging $\iota_i$ with $h_i \circ \iota$ replaces $H_i$ with a conjugate subgroup, we can suppose without loss of generality that $[\iota_i]$ intersects the interior of $N_X$. 
For each Enriques involution $\enr \in \aut(X)$ there exists a unique $i \in \{1,\ldots,r\}$ such that there exists $g\in G$ with $S_X^{\enr=1} = [\iota_i]^g$. 
Moreover, by Lemma~\ref{lem:w_acts_trivially}, we can suppose that $g \in \OG(S_X,N_X)$. 
We map such an $\enr$ to the double coset $H_iq(g)K \in H_i\backslash G / K$. 
This function is trivially well-defined and injective by Proposition~\ref{prop:main-Davide}.

To show surjectivity, take $i \in \{1,\ldots,r\}$ and $H_i\xi K \in H_i\backslash G / K$, with $\xi \in G$. By Proposition~\ref{prop:surjectivity}, $\xi = q(g)$ for some $g \in \OG(S_X,N_X)$. As $[\iota_i]^g$ also intersects the interior of $N_X$, by Proposition~\ref{prop:existence-ohashi} there is an Enriques involution $\enr \in \aut(X)$ which maps to $H_i\xi K$. This concludes the proof.
\endproof


\begin{corollary}
The number of Enriques involutions of a singular K3 surface~$X$ up to conjugation in $\aut(X)$ only depends on the genus of the transcendental lattice~$T_X$.
\end{corollary}
\proof The lattice $S_X$ is unique in its genus by Lemma~\ref{lem:surjectivity}, so it is completely determined by the genus of $T_X$. The subgroup $\OG(q(S_X),\omega_X)$ is also determined by the genus of $T_X$ when $X$ is singular (see Remark~\ref{rmk:OT-omega-singular-K3}). The subgroups $\OG(q(S_X),[\iota])$ for $\iota \in I_X$ only depend on $S_X$, so in turn they depend only on the genus of $T_X$.
\endproof

\begin{remark} \label{rmk:TXgenus}
Schütt~\cite{Schutt2007} described a relation of two singular K3 surfaces whose transcendental lattices are in the same genus.
See also~\cite{ShimadaReduction}.
\end{remark}

\subsection{Table \ref{tab:enriques_quotients} and Table \ref{tab:lattices}} \label{sec:tables}

Table~\ref{tab:enriques_quotients} contains the list of all singular K3 surfaces $X$ of discriminant $d$ with $d \leq 36$, given by their respective transcendental lattices $T_X$, together with the list of the Enriques involutions that they admit, up to conjugation in $\aut(X)$.

Table~\ref{tab:lattices} contains the list of the Gram matrices of all lattices appearing in Table~\ref{tab:enriques_quotients} which are not standard ADE lattices. 
The name $N^{\rho_2,\rho_4}_{r,d}$ (resp. $M^{\rho_2,\rho_4}_{r,d}$) denotes a positive definite even (resp. odd) lattice of rank $r$, determinant $d$, with $\rho_2$ $2$-vectors and $\rho_4$ $4$-vectors ($\rho_4$ omitted if not needed to distinguish two lattices). 
We use the following shorthand notation:
\begin{equation} \label{eq:shorthand}
[a_{11},a_{12},a_{22},\ldots,a_{1n},\ldots,a_{nn}] := 
\begin{bmatrix}
a_{11} & a_{12} & \ldots & a_{1n}\\
a_{12} & a_{22} & &\\
\vdots & & \ddots & \vdots \\
a_{1n} & & \cdots & a_{nn}
\end{bmatrix}. 
\end{equation}

\begin{center}
\footnotesize

\begin{longtable}{ccccccc}
  \caption{Enriques involutions up to conjugation of singular K3 surfaces of discriminant  $d \leq 36$.}
  \label{tab:enriques_quotients} \\
  
    \toprule
    $d$ & $T_X$ & $|\mathrm{Enr}|$ & $q(N)$ & $N$ & $|I_X(N)|$ & $\enr$ \\
    \midrule
  \endfirsthead

    \multicolumn{7}{c}%
    {\tablename\ \thetable{} -- continued from previous page} \\
    \midrule
    $d$ & $T$ & $|\mathrm{Enr}|$ & $q(N)$ & $N$ & $|I_X(N)|$ & $\enr$ \\
    \midrule
  \endhead

    \multicolumn{7}{c}{Continued on next page} \\ \midrule
  \endfoot

    \bottomrule
  \endlastfoot
  
    $3$ & $[2,1,2]$ & $0$ & -- & -- & -- \\
    \midrule
    $4$ & $[2,0,2]$ & $0$ & -- & -- & -- \\
    \midrule
    $7$ & $[2,1,4]$  & $2$ & $u_1^{\oplus 5} \oplus \anglefrac{2}{7}$ & $N^{144}_{10,7}(2)$ & $1$ & No.\,1 \\
               &               &&& $N^{242}_{10,7}(2)$ & $1$ & No.\,2 \\
    \midrule
    $8$ & $[2,0,4]$ & $0$ & -- & -- & -- \\
    \midrule
    $11$ & $[2,1,6]$ & $0$  & -- & -- & --\\
    \midrule
    $12$ & $[2,0,6]$ & $1$ & $u_1^{\oplus 4} \oplus \anglefrac{1}{2} \oplus \anglefrac{1}{6}$ & $M^{144}_{10,3}(2)$ & $1$ & No.\,3 \\
    \midrule
    $12$ & $[4,2,4]$ & $3$ & $u_1^{\oplus 4} \oplus v_1 \oplus \anglefrac{4}{3}$ & $N^{246}_{10,3}(2)$ & $3\times 1$ & Nos.\,\,4--6 \\
    \midrule
    $15$ & $[2,1,8]$ & $5$ & $u_1^{\oplus 5} \oplus \anglefrac{2}{15}$ & $N^{90}_{10,15}(2)$ & $1$ & No.\,$7$ \\
               &&&& $N^{132}_{10,15}(2)$ & $1$ & No.\,$11$ \\
               &&&& $N^{144}_{10,15}(2)$ & $2$ & Nos.\,$8,\,9$ \\
               &&&& $N^{240}_{10,15}(2)$ & $1$ & No.\,$10$ \\
    \midrule
    $15$ & $[4,1,4]$ & $4$ & $u_1^{\oplus 5} \oplus \anglefrac{4}{15}$ & $N^{92}_{10,15}(2)$  & $1$ \\
                &&&& $N^{112}_{10,15}(2)$ & $1$ \\
                &&&& $N^{242}_{10,15}(2)$ & $2$ \\
    \midrule
    $16$ & $[2,0,8]$ & $0$ & -- & -- & -- \\
    \midrule
    $16$ & $[4,0,4]$ & $9$ & $u_1^{\oplus 4} \oplus \anglefrac{1}{4} \oplus \anglefrac{1}{4}$ & $D_{10}(2)$ & $3 \times 1$ & Nos.\,12--14 \\
                &&&& $N^{244}_{10,4}(2)$   & $5 \times 1$ & Nos.\,15--18,\,20 \\
    \cmidrule{4-7}
                &&& $u_1^{\oplus 3} \oplus \anglefrac{1}{4} \oplus \anglefrac{1}{4}$ & $N^{0,308}_{10,1024}$ & $1$ & No.\,$19$ \\
    \midrule
    $19$ & $[2,1,10]$ & $0$ & -- & -- & -- \\
    \midrule
    $20$ & $[2,0,10]$ & $1$ & $u_1^{\oplus 4} \oplus \anglefrac{1}{2} \oplus \anglefrac{1}{10}$ & $M^{132}_{10,5}(2)$ & $1$ \\
    \midrule
    $20$ & $[4,2,6]$ & $2$ & $u_1^{\oplus 4} \oplus \anglefrac{3}{2} \oplus \anglefrac{3}{10}$ & $M^{92}_{10,5}(2)$ & $1$ & No.\,$21$ \\
                &&&& $M^{242}_{10,5}(2)$ & $1$ & No.\,$22$ \\
    \midrule
    $23$ & \multirow{2}{*}{\begin{tabular}{c} $[2,1,12],$ \\ $[4,\pm 1,6]$ \end{tabular}}  & $7$ & $u_1^{\oplus 5} \oplus \anglefrac{2}{23}$ & $N^{74}_{10,23}(2)$ & $1$ \\
                &&&& $N^{84}_{10,23}(2)$ & $1$ \\
                &&&& $N^{112}_{10,23}(2)$ & $1$ \\
                &&&& $N^{132}_{10,23}(2)$ & $1$ \\
                &&&& $N^{144}_{10,23}(2)$ & $1$ \\
                &&&& $N^{240}_{10,23}(2)$ & $1$ \\
                &&&& $N^{242}_{10,23}(2)$ & $1$ \\
    \midrule
    $24$ & $[2,0,12]$ & $1$ & $u_1^{\oplus 4} \oplus \anglefrac{1}{2} \oplus \anglefrac{1}{12}$ & $M^{90}_{10,6}(2)$ & $1$ & No.\,$23$ \\
    \midrule
    $24$ & $[4,0,6]$ & $1$ & $u_1^{\oplus 4} \oplus \anglefrac{3}{2} \oplus \anglefrac{11}{12}$ & $M^{242}_{10,6}(2)$ & $1$ \\
    \midrule
    $27$ & $[2,1,14]$ & $0$ & -- & -- & -- \\
    \midrule
    $27$ & $[6,3,6]$ & $0$ & -- & -- & -- \\
    \midrule
    $28$ & $[2,0,14]$ & $1$ & $u_1^{\oplus 4} \oplus \anglefrac{1}{2} \oplus \anglefrac{1}{14}$ & $M^{112}_{10,7}(2)$ & $1$ \\
    \midrule
    \newpage
    $28$ & $[4,2,8]$ & $24$ & $u_1^{\oplus 5} \oplus \anglefrac{2}{7}$ & $N^{144}_{10,7}(2)$ & $3\times 1 + 4\times 2$ \\
                &&&& $N^{242}_{10,7}(2)$ & $4\times 1+4\times 2$ \\
    \cmidrule{4-6}
                &&& $u_1^{\oplus 4} \oplus \anglefrac{2}{7}$ & $N^{0,274}_{10,1792}$ & $1$ \\
    \midrule
    $31$   & \multirow{2}{*}{\begin{tabular}{c} $[2,1,16],$ \\ $[4,\pm 1,8]$ \end{tabular}} & $9$ & $u_1^{\oplus 5} \oplus \anglefrac{2}{23}$ & $N^{60}_{10,31}(2)$ & $1$ \\
                &&&& $N^{72}_{10,31}(2)$ & $1$ \\
                &&&& $N^{86}_{10,31}(2)$ & $1$ \\
                &&&& $N^{90}_{10,31}(2)$ & $1$ \\
                &&&& $N^{112}_{10,31}(2)$ & $1$ \\
                &&&& $N^{128}_{10,31}(2)$ & $1$ \\
                &&&& $N^{144}_{10,31}(2)$ & $1$ \\
                &&&& $N^{240}_{10,31}(2)$ & $1$ \\
                &&&& $N^{242}_{10,31}(2)$ & $1$ \\ 
    \midrule
    $32$ & $[2,0,16]$ & $1$ & $u_1^{\oplus 4} \oplus \anglefrac{1}{2} \oplus \anglefrac{1}{16}$ & $M^{84}_{10,8}(2)$ & $1$ \\
    \midrule
    $32$ & $[4,0,8]$ & $33$ & $u_1^{\oplus 4} \oplus \anglefrac{1}{4} \oplus \anglefrac{1}{8}$ & $N^{138}_{10,8}(2)$ & $2\times 1 + 4 \times 2$ \\
                &&&& $N^{146}_{10,8}(2)$     & $3 \times 1 + 2 \times 2$ \\
                &&&& $N^{242}_{10,8}(2)$     & $3 \times 1 + 5 \times 2$ \\
    \cmidrule{4-6}
                &&& $u_1^{\oplus 3} \oplus \anglefrac{1}{4} \oplus \anglefrac{1}{8}$& $N^{0,210}_{10,2048}$  & $1$ \\
                &&&& $N^{0,250}_{10,2048}$  & $1$ \\
                &&&& $N^{0,274}_{10,2048}$  & $1$ \\
    \midrule
    $32$ & $[6,2,6]$ & $3$ & $u_1^{\oplus 4} \oplus \anglefrac{3}{2} \oplus \anglefrac{3}{16}$
    & $M^{112}_{10,8}(2)$ & $1$ \\
               &&&& $M^{144}_{10,8}(2)$ & $1$ \\
               &&&& $M^{240}_{10,8}(2)$ & $1$ \\
    \midrule
    $35$ & $[2,1,18]$ & $0$ & -- & -- & -- \\
    \midrule
    $35$ & $[6,1,6]$ & $0$ & -- & -- & -- \\
    \midrule
    $36$ & $[2,0,18]$ & $3$  & $u_1^{\oplus 4} \oplus \anglefrac{1}{2} \oplus \anglefrac{1}{18}$ & $M^{74}_{10,9}(2)$ & $1$\\
               &&&& $M^{90}_{10,9}(2)$ & $1$ \\
               &&&& $M^{128}_{10,9}(2)$ & $1$ \\
    \midrule
    $36$ & $[4,2,10]$ & $2$ & $u_1^{\oplus 4} \oplus \anglefrac{1}{2} \oplus \anglefrac{5}{18}$ & $M^{80}_{10,9}(2)$ & $1$ \\
               &&&& $M^{242}_{10,9}(2)$ & $1$ \\
    \midrule
    $36$ & $[6,0,6]$ & $3$ & $u_1^{\oplus 4} \oplus \anglefrac{1}{6} \oplus \anglefrac{1}{6}$ & $M^{60}_{10,9}(2)$ & $1$ & No.\,$24$ \\
               &&&& $M^{132}_{10,9}(2)$ & $1$ & No.\,$26$ \\
               &&&& $M^{240}_{10,9}(2)$ & $1$ & No.\,$25$ 
  \end{longtable}

\end{center}

\afterpage{\begin{landscape}
\begin{center}
\footnotesize
\begin{longtable}{cl}
  \caption{Lattices appearing in Table \ref{tab:enriques_quotients}.}
  \label{tab:lattices} \\
  
    \toprule
    Name & Gram matrix \\
    \midrule
  \endfirsthead

    \multicolumn{2}{c}%
    {\tablename\ \thetable{} -- continued from previous page} \\
    \midrule
    Name & Gram matrix \\
    \midrule
  \endhead

    \midrule
    \multicolumn{2}{c}{Continued on next page} \\ \midrule
  \endfoot

    \bottomrule
  \endlastfoot
  
  $M^{144}_{10,3}$ & $[ 2, 1, 2, 1, 0, 2, 0, 0, 1, 3, 1, 0, 0, 0, 2, 0, 0, 0, 0, 1, 2, 0, 0, 0, 0, 0, 1, 2, 0, 0, 0, 0, 0, 0, 1, 2, 0, 0, 0, 0, 0, 0, 0, 1, 2, 0, 0, 0, 0, 0, 0, 0, 0, 1, 2 ]$ \\
  $N^{246}_{10,3}$ & $[2,1,2] \oplus E_{8}$ \\
  $N^{244}_{10,4}$ & $[2] \oplus [2] \oplus E_{8}$ \\
  $M^{92}_{10,5}$ & $[ 2, 1, 2, 0, 1, 2, 0, 0, 1, 2, 0, 0, 0, 1, 2, 1, 0, 0, 0, 0, 2, 0, 0, 0, 0, 0, 1, 2, 0, 0, 0, 0, 0, 0, 1, 2, 0, 0, 0, 0, 0, 0, 0, 1, 2, 1, 0, 0, 0, 0, 0, 0, 0, 0, 3 ]$ \\
  $M^{132}_{10,5}$ & $[ 2, 1, 2, 1, 0, 2, 0, 0, 1, 2, 1, 0, 0, 0, 2, 0, 0, 0, 0, 1, 2, 0, 0, 0, 0, 0, 1, 2, 0, 0, 0, 0, 0, 0, 1, 3, 0, 0, 0, 0, 0, 0, 0, 1, 2, 0, 0, 0, 0, 0, 0, 0, 0, 1, 2 ]$ \\
  $M^{242}_{10,5}$ & $[2,1,3] \oplus E_{8}$ \\
  $M^{90}_{10,6}$ & $[ 2, 1, 2, 0, 1, 2, 0, 0, 1, 2, 0, 0, 0, 1, 2, 0, 0, 0, 0, 1, 2, 1, 0, 0, 0, 0, 0, 3, 1, 0, 0, 0, 0, 0, 0, 2, 0, 0, 0, 0, 0, 0, 0, 1, 2, 0, 0, 0, 0, 0, 0, 0, 0, 1, 2 ]$ \\
  $M^{242}_{10,6}$ & $[2] \oplus [3] \oplus E_{8}$ \\
  $M^{112}_{10,7}$ & $[ 2, 1, 2, 1, 0, 2, 0, 0, 1, 3, 1, 0, 0, 0, 2, 0, 0, 0, 0, 1, 2, 0, 0, 0, 0, 0, 1, 2, 0, 0, 0, 0, 0, 0, 1, 2, 0, 0, 0, 0, 0, 0, 0, 1, 2, 0, 0, 0, 0, 0, 0, 0, 0, 1, 3 ]$ \\
  $N^{144}_{10,7}$ & $[ 2, 1, 2, 1, 0, 2, 0, 0, 1, 4, 1, 0, 0, 0, 2, 0, 0, 0, 0, 1, 2, 0, 0, 0, 0, 0, 1, 2, 0, 0, 0, 0, 0, 0, 1, 2, 0, 0, 0, 0, 0, 0, 0, 1, 2, 0, 0, 0, 0, 0, 0, 0, 0, 1, 2 ]$ \\
  $N^{242}_{10,7}$ & $[2,1,4] \oplus E_{8}$ \\
  $M^{84}_{10,8}$ & $[ 2, 1, 2, 0, 1, 2, 1, 0, 0, 2, 0, 0, 0, 1, 2, 0, 0, 0, 0, 1, 3, 0, 0, 0, 0, 0, 1, 2, 0, 0, 0, 0, 0, 0, 1, 2, 0, 0, 0, 0, 0, 0, 1, 0, 2, 1, 0, 0, 0, 0, 0, 0, 0, 0, 2 ]$ \\
  $M^{112}_{10,8}$ & $[ 2, 1, 2, 0, 1, 2, 0, 0, 1, 2, 0, 0, 0, 1, 2, 0, 0, 0, 0, 1, 2, 1, 0, 0, 0, 0, 0, 2, 1, 0, 0, 0, 0, 0, 0, 2, 0, 0, 0, 0, 0, 0, 0, 1, 3, 0, 0, 0, 0, 0, 0, 0, 0, 1, 3 ]$ \\
  $N^{138}_{10,8}$ & $A_{3} \oplus E_{7}$  \\
  $M^{144}_{10,8}$ & $[ 2, 1, 2, 1, 0, 2, 1, 0, 0, 2, 0, 0, 0, 1, 2, 0, 0, 0, 0, 1, 2, 0, 0, 0, 0, 0, 1, 2, 0, 0, 0, 0, 0, 0, 1, 2, 0, 0, 0, 0, 0, 0, 0, 1, 2, 0, 0, 0, 0, 0, 0, 0, 0, 1, 3 ]$ \\
  $N^{146}_{10,8}$ & $[2] \oplus D_{9}$ \\
  $M^{240}_{10,8}$ & $[3,1,3] \oplus E_{8}$ \\
  $N^{242}_{10,8}$ & $[2] \oplus [4] \oplus E_{8}$ \\
  $M^{60}_{10,9}$ & $[ 2, 1, 2, 1, 0, 3, 0, 0, 1, 2, 0, 0, 0, 1, 2, 0, 0, 0, 0, 1, 2, 0, 0, 0, 1, 0, 0, 2, 0, 0, 0, 0, 0, 0, 1, 2, 1, 0, 0, 0, 0, 0, 0, 0, 2, 0, 0, 0, 0, 0, 0, 0, 0, 1, 2 ]$ \\
  $M^{74}_{10,9}$ & $[ 2, 1, 3, 0, 1, 2, 1, 0, 0, 2, 0, 0, 0, 1, 2, 0, 0, 0, 0, 1, 2, 0, 0, 0, 0, 0, 1, 2, 0, 0, 0, 0, 0, 0, 1, 2, 1, 0, 0, 0, 0, 0, 0, 0, 2, 0, 0, 0, 0, 0, 0, 0, 0, 1, 2 ]$ \\
  $M^{80}_{10,9}$ & $[ 3, 1, 2, 1, 0, 2, 0, 0, 1, 2, 0, 0, 0, 1, 2, 0, 0, 0, 0, 1, 2, 0, 0, 0, 0, 0, 1, 2, 1, 0, 0, 0, 0, 0, 0, 2, 0, 0, 0, 0, 0, 0, 0, 1, 2, 0, 0, 0, 0, 1, 0, 0, 0, 0, 2 ]$ \\
  $M^{90}_{10,9}$ & $[ 2, 1, 2, 0, 1, 2, 1, 0, 0, 2, 0, 0, 0, 1, 2, 0, 0, 0, 0, 1, 2, 0, 0, 0, 0,  0, 1, 2, 0, 0, 0, 0, 0, 0, 1, 2, 0, 0, 0, 0, 0, 0, 0, 1, 2, 1, 0, 0, 0, 0, 0, 0, 0, 0, 3 ]$ \\
  $M^{128}_{10,9}$ & $[ 2, 1, 2, 0, 1, 2, 0, 0, 1, 2, 0, 0, 0, 1, 3, 0, 0, 0, 0, 1, 2, 0, 0, 0, 0,  0, 1, 3, 1, 0, 0, 0, 0, 0, 0, 2, 0, 0, 0, 0, 0, 0, 0, 1, 2, 1, 0, 0, 0, 0, 0, 0, 0, 0, 2 ]$ \\
  $M^{132}_{10,9}$ & $[2,1,2] \oplus [ 2, 1, 2, 0, 1, 2, 0, 0, 1, 2, 0, 0, 0, 1, 3, 1, 0, 0, 0, 0, 2, 0, 0, 0, 0, 0, 1, 2, 1, 0, 0, 0, 0, 0, 0, 2 ]$ \\
  $M^{240}_{10,9}$ & $[3] \oplus [3] \oplus E_{8}$ \\
  $M^{242}_{10,9}$ & $[2,1,5] \oplus E_{8}$ \\
  
  
  $N^{90}_{10,15}$ & $[ 2, 1, 4, 1, 0, 2, 0, 0, 1, 2, 0, 0, 0, 1, 2, 0, 0, 0, 0, 1, 2, 1, 0, 0, 0,  0, 0, 2, 0, 0, 0, 0, 0, 0, 1, 2, 0, 0, 0, 0, 0, 0, 0, 1, 2, 0, 0, 0, 0, 0, 0, 0, 0, 1, 2 ]$ \\
  $N^{92}_{10,15}$ & $A_4 \oplus E_6$ \\
  $N^{112}_{10,15}$ & $[ 2, 1, 2, 0, 1, 2, 0, 0, 1, 2, 0, 0, 0, 1, 2, 0, 0, 0, 0, 1, 2, 0, 0, 0, 0, 0, 1, 2, 0, 0, 0, 0, 0, 1, 0, 4, -1, 0, 0, 0, 0, 0, 1, 0, 2, 1, 0, 0, 0, 0, 0, 0, 0, 0, 4 ]$ \\
  $N^{132}_{10,15}$ & $[2,1,2] \oplus [ 2, 1, 2, 1, 0, 2, 0, 0, 1, 2, 1, 0, 0, 0, 2, 0, 0, 0, 0, 1, 2, 0, 0, 0, 0, 0, 1, 2, 0, 0, 0, 0, 0, 0, 1, 4 ]$ \\
  $N^{144}_{10,15}$ & $[ 2, 1, 2, 1, 0, 2, 0, 0, 1, 6, 1, 0, 0, 0, 2, 0, 0, 0, 0, 1, 2, 0, 0, 0, 0, 0, 1, 2, 0, 0, 0, 0, 0, 0, 1, 2, 0, 0, 0, 0, 0, 0, 0, 1, 2, 0, 0, 0, 0, 0, 0, 0, 0, 1, 2 ]$ \\
  $N^{240}_{10,15}$ & $[4,1,4] \oplus E_{8}$ \\
  $N^{242}_{10,15}$ & $[2,1,8] \oplus E_{8}$ \\
  $N^{74}_{10,23}$ & $[ 2, 1, 2, 0, 1, 2, 0, 0, 1, 2, 0, 0, 0, 1, 2, 1, 0, 0, 0, 0, 2, 0, 0, 0, 0, 0, 1, 2, 0, 0, 0, 0, 0, 0, 1, 2, 1, 0, 0, 0, 0, 0, 0, 0, 4, 0, 0, 0, 0, 0, 0, 0, 0, 1, 2 ]$ \\
  $N^{84}_{10,23}$ & $[ 2, 1, 2, 0, 1, 2, 1, 0, 0, 2, 1, 0, 0, 0, 2, 0, 0, 0, 0, 1, 2, 0, 0, 0, 0, 0, 1, 4, 0, 0, 0, 0, 0, 0, 1, 2, 0, 0, 0, 0, 0, 0, 0, 1, 2, 0, 0, 0, 0, 0, 0, 0, 0, 1, 2 ]$ \\
  $N^{112}_{10,23}$ & $[ 2, 1, 2, 1, 0, 2, 0, 0, 1, 4, 1, 0, 0, 0, 2, 0, 0, 0, 0, 1, 2, 0, 0, 0, 0, 0, 1, 2, 0, 0, 0, 0, 0, 0, 1, 2, 0, 0, 0, 0, 0, 0, 0, 1, 2, 0, 0, 0, 0, 0, 0, 0, 0, 1, 4 ]$ \\
  $N^{132}_{10,23}$ & $[ 2, 1, 2, 1, 0, 2, 0, 0, 1, 2, 1, 0, 0, 0, 2, 0, 0, 0, 0, 1, 2, 0, 0, 0, 0, 0, 1, 2, 0, 0, 0, 0, 0, 0, 1, 6, 0, 0, 0, 0, 0, 0, 0, 1, 2, 0, 0, 0, 0, 0, 0, 0, 0, 1, 2 ]$ \\
  $N^{144}_{10,23}$ & $[ 2, 1, 2, 1, 0, 2, 0, 0, 1, 8, 1, 0, 0, 0, 2, 0, 0, 0, 0, 1, 2, 0, 0, 0, 0, 0, 1, 2, 0, 0, 0, 0, 0, 0, 1, 2, 0, 0, 0, 0, 0, 0, 0, 1, 2, 0, 0, 0, 0, 0, 0, 0, 0, 1, 2 ]$ \\
  $N^{240}_{10,23}$ & $[4,1,6] \oplus E_{8}$ \\ 
  $N^{242}_{10,23}$ & $[2,1,12] \oplus E_{8}$ \\
  $N^{60}_{10,31}$ & $[ 2, 1, 2, 1, 0, 2, 0, 0, 1, 2, 1, 0, 0, 0, 4, 0, 0, 0, 0, 1, 2, 0, 0, 0, 0, 0, 1, 2, 0, 0, 0, 0, 0, 0, 1, 2, 0, 0, 0, 0, 0, 0, 1, 0, 2, 0, 0, 0, 0, 0, 0, 0, 0, 1, 2 ]$ \\
  $N^{72}_{10,31}$ & $[ 2, 1, 4, 0, 2, 4, 0, 0, 1, 2, 0, 0, 0, 1, 2, 1, 0, 0, 0, 0, 2, 0, 0, 0, 0, 0, 1, 2, -1, 0, 0, 0, 0, 0, 0, 2, 0, 0, 0, 0, 0, 0, 0, 1, 2, 0, 0, 0, 0, 1, 0, 0, 0, 1, 2 ]$ \\
  $N^{86}_{10,31}$ & $[ 2, 1, 2, 0, 1, 2, 0, 0, 1, 2, 0, 0, 0, 1, 2, 1, 0, 0, 0, 0, 2, 1, 0, 0, 0, 0, 0, 2, 0, 0, 0, 0, 0, 0, 1, 4, 0, 0, 0, 0, 0, 0, 0, 2, 4, 0, 0, 0, 0, 0, 0, 0, 0, 1, 2 ]$ \\
  $N^{90}_{10,31}$ & $[ 2, 1, 2, 0, 1, 2, 0, 0, 1, 2, 0, 0, 0, 1, 2, 0, 0, 0, 0, 1, 2, 0, 0, 0, 0, 0, 1, 2, 0, 0, 0, 0, 0, 0, 1, 2, 0, 0, 0, 0, 0, 0, 0, 1, 2, 1, 0, 0, 0, 0, 0, 0, 0, 0, 4 ]$ \\
  $N^{112}_{10,31}$ & $[ 2, 1, 2, 0, 1, 2, 1, 0, 0, 4, 0, 0, 0, 1, 2, 0, 0, 0, 0, 1, 2, 0, 0, 0, 0, 0, 1, 2, -1, 0, 0, 0, 0, 0, 0, 2, 0, 0, 0, 0, 0, 1, 0, 1, 2, 0, 0, 0, 0, 1, 0, 0, 0, 0, 6 ]$ \\
  $N^{128}_{10,31}$ & $[ 2, 1, 2, 0, 1, 2, 1, 0, 0, 2, 0, 0, 0, 1, 2, 0, 0, 0, 0, 1, 2, 1, 0, 0, 0,  0, 0, 2, 0, 0, 0, 0, 0, 1, 0, 4, 0, 0, 0, 0, 0, 0, 0, 1, 4, 0, 0, 0, 0, 0, 0, 0, 0, 1, 2 ]$ \\
  $N^{144}_{10,31}$ & $[ 2, 1, 2, 0, 1, 2, 0, 0, 1, 2, 0, 0, 0, 1, 2, 0, 0, 0, 0, 1, 2, 0, 0, 0, 0, 0, 1, 2, 1, 0, 0, 0, 0, 0, 0, 2, 1, 0, 0, 0, 0, 0, 0, 0, 2, 0, 0, 0, 0, 0, 0, 0, 0, 1, 10 ]$ \\
  $N^{240}_{10,31}$ & $[4,1,8] \oplus E_{8}$ \\ 
  $N^{242}_{10,31}$ & $[2,1,16] \oplus E_{8}$ \\  

  $N^{0,308}_{10,1024}$ & $[ 4, 2, 4, 0, 2, 4, 0, 0, 1, 4, 0, 0, 0, 2, 4, 2, 0, 0, 0, 0, 4, 2, 0, 0, 0, 0, 0, 4, 0, 0, 0, 0, 0, 0, 2, 4, 0, 0, 0, 0, 0, 0, 0, 2, 4, 0, 0, 0, 0, 0, 0, 0, 0, 2, 4 ]$ \\
  $N^{0,274}_{10,1792}$ & $[ 4, 2, 4, 0, 2, 4, 0, 0, 1, 4, 0, 0, 0, 1, 4, 2, 0, 0, 0, 0, 4, 0, 0, 0, 0, 0, 2, 4, 0, 0, 0, 0, 0, 0, 2, 4, 0, 0, 0, 0, 0, 0, 0, 2, 4, 2, 0, 0, 0, 0, 0, 0, 0, 0, 4 ]$ \\
  $N^{0,210}_{10,2048}$ & $[ 4, 2, 4, 2, 0, 4, 2, 0, 0, 4, 0, 0, 0, 2, 4, 0, 0, 0, 0, 2, 4, 0, 0, 0, 0, 0, 2, 4, 0, 0, 0, 0, 0, 0, 2, 4, 0, 0, 0, 0, 0, 0, 0, 2, 4, 0, 0, 0, 1, 0, 0, 0, 0, 0, 4 ]$ \\
  $N^{0,250}_{10,2048}$ & $[ 4, 2, 4, 1, 0, 4, 0, 0, 1, 4, 0, 0, 0, 2, 4, 0, 0, 0, 0, 2, 4, 0, 0, 0, 0, 0, 2, 4, 2, 0, 0, 0, 0, 0, 0, 4, 0, 0, 0, 0, 0, 2, 0, 2, 4, 0, 0, 2, 0, 0, 0, 0, 0, 0, 4 ]$ \\ 
  $N^{0,274}_{10,2048}$ & $[4] \oplus [ 4, 2, 4, 0, 2, 4, 0, 0, 2, 4, 1, 0, 0, 0, 4, 0, 0, 0, 0, 2, 4, 0, 0, 0, 0, 0, 2, 4, 1, 0, 0, 0, -2, 0, 0, 4, 0, 0, 0, 0, 0, 0, 0, 2, 4 ]$ \\
\end{longtable}
\end{center}
\end{landscape}}

We will illustrate presently how these two tables were compiled.

The following theorem by Sertöz builds on work by Keum~\cite{Keum1990} and characterizes singular K3 surfaces without Enriques quotients. 
\begin{theorem}[Sertöz~\cite{Sertoz2005}; see also \cite{HS2012}] \label{thm:sertoez}
Let $X$ be a singular K3 surface of discriminant $d$. Then $X$ has no Enriques involution if and only if
$d \equiv 3\,(8)$ or $T_X \in \{[2,0,2],[2,0,4],[2,0,8]\}$.
\noproof
\end{theorem}
In all other cases, we determined the set of conjugacy classes of all Enriques involutions in $\aut(X)$ by means of Theorem~\ref{thm:davide-main}. The item $|\mathrm{Enr}|$ indicates the number of such conjugacy classes.

First of all, one must determine a complete set of representatives for the action of $\OG(S_X)$ on $I_X$. Given a positive definite even lattice $N$ of rank $10$ without $2$-vectors (see Theorem~\ref{thm:keum}), we put
\[
 I_X(N) := \shortset{\iota \in I_X}{[\iota]^\perp \cong N(-1)}.
\]
Clearly, the sets $I_X(N)$ form a partition of $I_X$ which respects the $\OG(S_X)$-action, so we reduce the problem to computing a complete set of representatives for the action of $\OG(S_X)$ on $I_X(N)$, for each $N$ such that $I_X(N) \neq \emptyset$. 

We find all such lattices in the following way. 
Using Proposition~\ref{prop:Nikulin1.15}, we list all possible finite quadratic forms $q$, such that $q \cong q(N)$.
For each form $q$, we determine all lattices $N$ in the genus $\fg(10,0,q)$ without $2$-vectors (see Algorithm~\ref{alg:lattices-in-genus}).
All possible finite quadratic forms $q = q(N)$ and orthogonal complements $N$ have been listed in Table~\ref{tab:enriques_quotients}.

Since $I_X(N) = I(S_X,L_{10}(2),N(-1))$ as defined in Section~\ref{sec:primitive-sublattices}, a complete set of representatives $\iota_{1},\ldots,\iota_r$ up to the action of $\OG(S_X)$ on $I_X(N)$ can be enumerated using Proposition \ref{prop:Nikulin1.5.1}.
For each $i \in \{1,\ldots,r\}$, the subgroup $H_i = \OG(q(S_X),[\iota_i])$ of $G = \OG(q(S_X))$ can be determined using Proposition~\ref{prop:OMS}. 
On the other hand, the subgroup $K = \OG(q(S_X),\omega_X)$ can be computed using Remark~\ref{rmk:OT-omega-singular-K3}.

\begin{remark}
In order to apply Proposition \ref{prop:Nikulin1.5.1}, it is worth mentioning that for $L = L_{10}(2)$ the natural homomorphism $\OG(L) \rightarrow \OG(q(L))$ is surjective and that, up to the action of $\OG(q(L))$, there are only two subgroups of $L^\vee/L$ of order $2$.

On the other hand, since $N$ is positive definite, we can compute $\OG(N)$ by the attribute \verb+automorphism_group+ of the class \verb+QuadraticForm+ in {\tt sage}; hence, we can compute its image in $\OG(q(N))$. 
Such a function has also been implemented for the class \verb+IntegralLattice+ by Brandhorst~\cite{Brandhorst}.
\end{remark}

The item $|I_X(N)|$ gives the cardinalities of the sets of double cosets $H_i\backslash G / K$. For instance, the entry ``$3\times 1 + 4\times 2$'' means that $r = 7$, $|H_i\backslash G / K| = 1$ for $i =1,2,3$ and $|H_i\backslash G / K| = 2$ for $i = 4,\ldots,7$. 
Note that the item $|\mathrm{Enr}|$ is the sum of the items $|I_X(N)|$ over the lattices $N$.

Finally, the item $\enr$ refers to those involutions studied in detail in Section~\ref{sec:enriques-involutions-of-certain-K3} and listed in Table~\ref{table:enrtable}.


%
\section{Automorphism groups of singular K3 surfaces}\label{sec:AutSingK3}
\subsection{Borcherds method}\label{subsec:Borcherds}
We explain Borcherds method~(\cite{Bor1},~\cite{Bor2})
to calculate $\aut(X)$ of a K3 surface $X$
and its action on $N_X$.
The details of the algorithms in the
computation below are explained in~\cite{ShimadaAlgo}.
Suppose that we have a primitive embedding
\[
\iotaX\colon S_X \inj L_{26}.
\]
We assume that $\iotaX$ maps $\PPP_X$ to the positive cone $\PPP_{26}$ of $L_{26}$,
and consider the decomposition of $\PPP_X$ by 
$\iotaX^*\RRR_{26}\sperp$-chambers, that is, by chambers induced by 
Conway chambers non-degenerate with respect to $\iotaX$.
Since $\iotaX$ maps $\RRR_X$ to $\RRR_{26}$,
every $\RRR_X\sperp$-chamber is a union of $\Xinduced$-chambers.
In particular, the nef chamber $N_X$ is a union of $\Xinduced$-chambers.
Since a Conway chamber is quasi-finite,
every $\Xinduced$-chamber is quasi-finite.
\par
The orthogonal complement $[\iotaX]\sperp$ of the image of $\iotaX$ 
is an even negative definite lattice.
The even unimodular overlattice $L_{26}$ of $S_X\oplus [\iotaX]\sperp$ 
induces an \emph{anti}-isometry $q(S_X)\cong -q([\iotaX]\sperp)$,
and hence an isomorphism $\OG(q(S_X))\cong\OG(q([\iotaX]\sperp))$.
We assume the following condition:
\[
\parbox{11cm}{\textrm{the image of $\OG(q(S_X), \omega_X)$ by 
the isomorphism $\OG(q(S_X))\cong\OG(q([\iotaX]\sperp))$ above
is contained in the image 
of the natural homomorphism $ \OG([\iotaX]\sperp) \to \OG(q([\iotaX]\sperp))$.}} \tag{A} \label{eq:condA}
\]
Since $\OG([\iotaX]\sperp)$ and $\OG(q(S_X), \omega_X)$ are finite, 
we can determine whether this condition is fulfilled or not.
Suppose that Condition~\eqref{eq:condA} is satisfied.
Then every isometry $g\in \OG(S_X, \omega_X) \cap \OG(S_X,\cP_X)$ extends to 
an isometry $\tilde{g}\in \OG(L_{26}, \PPP_{26})$,
which preserves the set of Conway chambers.
Therefore every isometry of $S_X$ satisfying the period condition 
preserves the set of $\Xinduced$-chambers.
\par
We also assume the following condition:
\[
\textrm{$[\iotaX]\sperp$ cannot be embedded into the negative definite Leech lattice.} \tag{B} \label{eq:condB}
\]
For example, if $[\iotaX]\sperp$ contains a $(-2)$-vector, then this condition is fulfilled.
Condition~\eqref{eq:condB} implies that each $\Xinduced$-chamber $D$ in $\PPP_X$ 
has only a finite number of walls (see~\cite{ShimadaAlgo}).
More precisely,
if $D$ is induced by a Conway chamber $C$,
then the set of vectors defining walls of $D$ can be calculated 
from the Weyl vector $\weyl_C$ corresponding to $C$
by Theorem~\ref{thm:Conway}.
By this finiteness, 
we can calculate,
for two $\Xinduced$-chambers $D$ and $D\sprime$,
the set of all isometries
$g\in \OG(S_X)$ such that $D^{g}=D\sprime$.
In particular,
the group
\[
\OG(S_X, D):=\shortset{g\in \OG(S_X)}{D^g=D}
\]
is finite, and can be calculated explicitly.
If $D\subset N_X$, then 
\[
\aut(X, D):=\OG(S_X, D)\cap \OG(S_X, \omega_X)
\]
is contained in $\aut(X)$, and can be calculated explicitly.
\begin{definition}
Let $D$ be an $\Xinduced$-chamber contained in $N_X$.
A wall $D\cap (v)\sperp$ of $D$ is called an \emph{outer wall} 
if it is defined by a $(-2)$-vector, that is, if 
there exists a rational number $\lambda$ such that $-2/\intf{v, v}=\lambda^2$
and $\lambda v\in S_X$.
Otherwise, we say that $D\cap (v)\sperp$ is an \emph{inner wall}.
\end{definition}
A wall $D\cap (v)\sperp$ is an outer wall if and only if $N_X\cap (v)\sperp$
is a wall of $N_X$.
The $\Xinduced$-chamber $D\sprime$adjacent to $D$ across a wall $D\cap (v)\sperp$ of $D$
is contained in $N_X$ if and only if $D\cap (v)\sperp$ is an inner wall.
\par
Let $D$ be an $\Xinduced$-chamber,
and let $\weyl_C$ be the Weyl vector corresponding to a Conway chamber $C$ inducing 
$D=\iota_X\inv(C)$.
Let $D\cap (v)$ be a wall of $D$,
and let $D\sprime$ be 
the $\Xinduced$-chamber adjacent to $D$ across $D\cap (v)\sperp$.
Then we can calculate the Weyl vector $\weyl_{C\sprime}$ 
corresponding to a Conway chamber $C\sprime$ inducing $D\sprime=\iota_X\inv(C\sprime)$
(see~\cite{ShimadaAlgo}), 
and hence we can calculate the set of walls of $D\sprime$,
which is again finite.
Therefore
we can determine whether there exists an isometry $g\in \OG(S_X, \omega_X)$
that maps $D$ to $D\sprime$. 
\begin{definition}\label{def:extraX}
Let $D\cap (v)\sperp$ be an inner wall of
an $\Xinduced$-chamber $D$ contained in $N_X$.
An isometry $g\in \OG(S_X, \omega_X)$
is said to be an \emph{extra automorphism} 
associated with $D\cap (v)\sperp$
if $g$ maps $D$ to the $\Xinduced$-chamber adjacent to 
$D$ across $D\cap (v)\sperp$.
\end{definition}
Let $g$ be an extra automorphism as above.
Since $g$ satisfies the period condition,
Condition~\eqref{eq:condA} implies that 
$g$ preserves the set of $\Xinduced$-chambers.
Moreover $g$ maps an interior point of $N_X$ to the interior of $N_X$,
and hence $g\in \aut(X)$.
We consider the following condition:
\[
\parbox{11cm}{There exists an $\Xinduced$-chamber $D_0$ contained in $N_X$
such that every \\
inner wall of $D_0$ has an extra automorphism.} \tag{IX} \label{eq:condIX}
\]
\begin{definition}\label{def:simple}
We say that an embedding $\iotaX$ satisfying Conditions~\eqref{eq:condA}, \eqref{eq:condB} and \eqref{eq:condIX} is 
of \emph{simple Borcherds type}.
\end{definition}
%
%
\begin{theorem}[\cite{ShimadaAlgo}]\label{thm:D0}
Suppose that $\iotaX$ is of simple Borcherds type.
\begin{enumerate}[{\rm (1)}]
\item For any point $v$ of $N_X$, 
there exists an automorphism $g$ of $X$ such that $v^g\in D_0$.
\item Let $o_1, \dots, o_m$ be the orbits of the action of $\aut(X, D_0)$
on the set of inner walls of $D_0$,
and, for $i=1, \dots, m$, 
let $g(o_i)$ be an extra automorphism associated with 
an inner wall $D_0\cap (v_i)\sperp$ belonging to $o_i$.
Then $\aut(X)$ is generated by 
$\aut(X, D_0)$
and the extra automorphisms $g(o_1), \dots, g(o_m)$.
\noproof
\end{enumerate}
\end{theorem}
\subsection{Application to certain singular K3 surfaces}\label{subsec:singK3}
We consider singular K3 surfaces with transcendental lattice 
$[a, b, c]$
in Table~\ref{table:Table1}.
These transcendental lattices are characterized 
among all even binary positive definite lattices by the following properties:
there exists a primitive embedding $\iotaX\colon S_X \inj L_{26}$ of 
simple Borcherds type such that 
the orthogonal complement $[\iotaX]\sperp$ is 
generated by $(-2)$-vectors.
In particular, Condition \eqref{eq:condB} is satisfied.
The column 
$\texttt{root}\,\texttt{type}$ in Table~\ref{table:Table1} 
indicates the $\ADE$-type of the standard fundamental root system 
of $[\iotaX]\sperp$.
For these cases, 
the natural homomorphism $\OG([\iotaX]\sperp)\to \OG(q([\iotaX]\sperp))$ is surjective
and hence Condition \eqref{eq:condA} is satisfied.
\begin{table}
    \caption{$T_X$ and $[\iotaX]\sperp$.}
    \label{table:Table1}
    \centering
    \begin{tabular}{ccccccccc}
    \toprule
    No. & $T_X$ & \texttt{root type} & $m_1$ & $m_2$ & $m_3$ & $m_4$ & $k_1$ & $k_2$ \\ 
    \midrule
    $1$ & $[2, 1, 2]$ & $E_{6}$ & $12$ & $6$ & $6$ & $3$ & $103680$ & $2$  \\ 
    $2$ & $[2, 0, 2]$ & $D_{6}$ & $8$ & $4$ & $4$ & $2$ & $46080$ & $2$  \\ 
    $3$ & $[2, 1, 4]$ & $A_{6}$ & $4$ & $2$ & $2$ & $1$ & $10080$ & $2$  \\ 
    $4$ & $[2, 0, 4]$ & $D_{5}+A_{1}$ & $4$ & $2$ & $2$ & $1$ & $7680$ & $2$  \\ 
    $5$ & $[2, 0, 6]$ & $A_{5}+A_{1}$ & $4$ & $2$ & $2$ & $1$ & $2880$ & $2$  \\ 
    $6$ & $[4, 2, 4]$ & $D_{4}+A_{2}$ & $12$ & $6$ & $1$ & $1$ & $13824$ & $12$ \\ 
    $7$ & $[2, 1, 8]$ & $A_{4}+A_{2}$ & $4$ & $2$ & $2$ & $1$ & $2880$ & $4$ \\ 
    $8$ & $[4, 0, 4]$ & $2A_{3}$ & $8$ & $4$ & $1$ & $1$ & $4608$ & $8$ \\ 
    $9$ & $[4, 2, 6]$ & $A_{4}+2A_{1}$ & $4$ & $2$ & $1$ & $1$ & $1920$ & $4$ \\
    $10$ & $[2, 0, 12]$ & $A_{3}+A_{2}+A_{1}$ & $4$ & $2$ & $2$ & $1$ & $1152$ & $4$ \\ 
    $11$ & $[6, 0, 6]$ & $2A_{2}+2A_{1}$ & $8$ & $4$ & $1$ & $1$ & $2304$ & $16$  \\
    \bottomrule
    \end{tabular}
\end{table}
%
The following data are also given in Table~\ref{table:Table1}.
\begin{itemize}
\item $m_1$ is the order of $\OG(T_X)$, $m_2$ is the order of $\OG(T_X, \omega_X)$,
$m_3$ is the order of the kernel $K$ of the homomorphism $\OG(T_X)\to \OG(q(T_X))$,
and $m_4$ is the order of $\OG(T_X, \omega_X)\cap K$.
Then $m_4$ is the order of the kernel of $\rho_X$ by Remark~\ref{rem:Kerrho},
and the order of $\OG(q(T_X), \omega_X)\cong \OG(q(S_X), \omega_X)$ is $m_2/m_4$.
\item $k_1$ is the order of $\OG([\iotaX]\sperp)$,
and $k_2$ is the order of $\OG(q(T_X)) \cong \OG(q(S_X)) \cong \OG(q([\iotaX]\sperp))$.
\end{itemize}
We have a Conway chamber $C_0$ that induces an $\Xinduced$-chamber $D_0$ contained in~$N_X$.
Let $\weyl\in L_{26}$ be the Weyl vector corresponding to $C_0$, 
and let $\weyl_S\in S_X\tensor\Q$ be the image of $\weyl$
by the orthogonal projection $\pr_S\colon L_{26}\tensor\Q\to S_X\tensor\Q$.
For each of the $11$ cases,
we can confirm that $\weyl_S$ belongs to the interior of $D_0$
and that $\weyl_S$ is invariant under the action of $\aut(X, D_0)$.
Let $o$ be an orbit of the action of $\aut(X, D_0)$
on the set of walls of $D_0$,
and let $D_0\cap (v)\sperp$ be a member of $o$.
We choose the defining vector $v$ of this wall 
in such a way that  $v$ is \emph{primitive} in $S_X\dual$.
Then $v$ is unique. The values 
$n:=\intf{v, v}$ and $a:=\intf{v, \weyl_S}$
are independent of the choice of the wall $D_0\cap (v)\sperp\in o$.
Suppose that the orbit $o$ consists of inner walls.
Then we can find an extra automorphism $g\in \aut(X)$ associated with $D_0\cap (v)\sperp$
by a direct calculation.
Hence $\iota_X$ is of simple Borcherds type.
The degree 
$d_g:=\intf{\weyl_S^g, \weyl_S}$
is also independent of the choice of $D_0\cap (v)\sperp$ and $g$.
Table~\ref{tab:D0data} contains the data of walls and extra automorphisms of $D_0$.
If $D_0\cap (v)\sperp$ is an inner wall,
the $(-2)$-vectors $r$ of $L_{26}$
such that $(r)\sperp$ passes through $\iota_X(D_0\cap (v)\sperp)\subset \PPP_{26}$
form a root system, whose
$\ADE$-type is also given below.
\begin{remark}
Almost all results in Table~\ref{tab:D0data} have already appeared in previous works.
See Vinberg~\cite{Vinberg1983} for Nos.\,1 and 2 of Table~\ref{table:Table1},
Ujikawa~\cite{Ujikawa2013} for No.\,3,
Keum and Kondo~\cite{KeumKondo2001} for Nos.\,6 and 8, 
\cite{ShimadaAlgo} for Nos.\,4, 5 and 6, 
\cite{ShimadaSch} for Nos.\,7, 9 and 11.
\end{remark}
\begin{remark}
In Table~\ref{tab:D0data},  the order of $\aut_0:=\aut(X, D_0)$ is given.
 We give the list of all elements of the finite group 
 $\aut(X, D_0)$ in~\cite{thecompdata},
 and hence we can determine its group structure.
 For example,
 for the case $T_X=[4,0,4]$,
 we see that  $\aut(X, D_0)$ is isomorphic to $(C_2^5) : S_5$.
\end{remark}
%


\begin{center}
\small

\begin{longtable}{cccrcccccc}
   \caption{Walls and extra automorphisms of $D_0$.}
  \label{tab:D0data} \\
  
    \toprule
    $T_X$ & $|\aut_0|$ & $\intf{\weyl_S, \weyl_S}$ & No. & $|o|$ & & $n$ & $a$ & $d_g$ & root type \\
    \midrule
  \endfirsthead

    \multicolumn{10}{c}%
    {\tablename\ \thetable{} -- continued from previous page} \\
    \midrule
    $T_X$ & $|\aut_0|$ & $\intf{\weyl_S, \weyl_S}$ & No. & $|o|$ & & $n$ & $a$ & $d_g$ & root type \\
    \midrule
  \endhead

    \multicolumn{10}{c}{Continued on next page} \\ \midrule
  \endfoot

    \bottomrule
  \endlastfoot
  
    $[2,1,2]$ & $72$ & $78$ & & $6$ & outer & $-2$ & $1$ & \\
    &&& & $18$ & outer & $-2$ & $1$ & \\
    \cmidrule{4-10}
    &&& $1$ & $12$ & inner & $-2/3$ & $9$ & $321$ & $E_{7}$ \\
    \midrule
    $[2,0,2]$ & $120$ & $55$ & & $10$ & outer & $-2$ & $1$ \\ 
    &&& & $15$ & outer & $-2$ & $1$ & \\ 
    &&& & $20$ & outer & $-1/2$ & $17/2$ & \\ 
    \cmidrule{4-10}
    &&& $1$ & $5$ & inner & $-1$ & $6$ & $127$ & $D_{7}$\\
    \midrule
    $[2,1,4]$ & $336$ & $28$ & & 28 & outer & $-2$ & $1$ \\
    \cmidrule{4-10}
    &&& $1$ & $14$ & inner & $-8/7$ & $4$ & $56$ & $A_{7}$ \\ 
    &&& $2$ & $28$ & inner & $-4/7$ & $6$ & $154$ & $D_{7}$ \\
    &&& $3$ & $56$ & inner & $-2/7$ & $7$ & $371$ & $E_{7}$ \\ 
    \midrule
    $[2,0,4]$ & $48$ & $61/2$ & & $6$ & outer & $-2$ & $1$ & \\ 
    &&& & $8$ & outer & $-2$ & $1$ & \\ 
    &&& & $12$ & outer & $-2$ & $1$ & \\ 
    &&& & $2$ & outer & $-1/2$ & $11/2$ & \\ 
    \cmidrule{4-10}
    &&& $1$ & $3$ & inner & $-3/2$ & $3/2$ & $67/2$ & $A_{2}+D_{5}$ \\ 
    &&& $2$ & $4$ & inner & $-1$ & $5$ & $161/2$ & $A_{1}+D_{6}$ \\ 
    &&& $3$ & $6$ & inner & $-1$ & $5$ & $161/2$ & $A_{1}+D_{6}$ \\ 
    &&& $4$ & $8$ & inner & $-3/4$ & $6$ & $253/2$ & $A_{1}+E_{6}$ \\ 
    &&& $5$ & $24$ & inner & $-3/4$ & $6$ & $253/2$ & $A_{1}+E_{6}$ \\ 
    &&& $6$ & $8$ & inner & $-1/4$ & $13/2$ & $737/2$ & $E_{7}$ \\
    \midrule
    $[2,0,6]$ & $144$ & $18$ & & $12$ & outer & $-2$ & $1$ & \\ 
    &&& & $18$ & outer & $-2$ & $1$ & \\ 
    &&& & $12$ & outer & $-1/2$ & $11/2$ & \\ 
    &&& & $36$ & outer & $-1/2$ & $11/2$ & \\ 
    \cmidrule{4-10}
    &&& $1$ & $4$ & inner & $-3/2$ & $3/2$ & $21$ & $A_{2}+A_{5}$ \\ 
    &&& $2$ & $24$ & inner & $-7/6$ & $7/2$ & $39$ & $A_{1}+A_{6}$ \\ 
    &&& $3$ & $6$ & inner & $-2/3$ & $4$ & $66$ & $A_{7}$ \\
    &&& $4$ & $24$ & inner & $-2/3$ & $5$ & $93$ & $A_{1}+D_{6}$ \\ 
    &&& $5$ & $36$ & inner & $-2/3$ & $5$ & $93$ & $A_{1}+D_{6}$ \\ 
    &&& $6$ & $24$ & inner & $-1/6$ & $11/2$ & $381$ & $E_{7}$ \\ 
    \midrule 
    $[4,2,4]$ & $1152$
    & $16$ & & $32$ & outer & $-2$ & $1$ & \\ 
    \cmidrule{4-10}
    &&& $1$ & $8$ & inner & $-4/3$ & $2$ & $22$ & $A_{3}+D_{4}$ \\ 
    &&& $2$ & $72$ & inner & $-1$ & $4$ & $48$ & $A_{2}+D_{5}$ \\ 
    &&& $3$ & $96$ & inner & $-1/3$ & $5$ & $166$ & $D_{7}$ \\ 
    \midrule
    $[2,1,8]$ & $720$ & $12$ & & $36$ & outer & $-2$ & $1$ & \\ \cmidrule{4-10}
    &&& $1$ & $12$ & inner & $-4/3$ & $2$ & $18$ & $A_{3}+A_{4}$ \\ 
    &&& $2$ & $40$ & inner & $-6/5$ & $3$ & $27$ & $A_{2}+A_{5}$ \\ 
    &&& $3$ & $90$ & inner & $-4/5$ & $4$ & $52$ & $A_{2}+D_{5}$ \\ 
    &&& $4, 5$ & $30$ & inner & $-8/15$ & $4$ & $72$ & $A_{7}$ \\ 
    &&& $6, 7$ & $120$ & inner & $-2/15$ & $5$ & $387$ & $E_{7}$ \\ 
    %
    \midrule
    $[4,0,4]$ & $3840$ & $10$ & & $40$ & outer & $-2$ & $1$ & \\ 
    \cmidrule{4-10}
    &&& $1$ & $64$ & inner & $-5/4$ & $5/2$ & $20$ & $A_{3}+A_{4}$ \\ 
    &&& $2$ & $40$ & inner & $-1$ & $3$ & $28$ & $A_{3}+D_{4}$ \\ 
    &&& $3$ & $160$ & inner & $-1/2$ & $4$ & $74$ & $A_{7}$ \\ 
    &&& $4$ & $320$ & inner & $-1/4$ & $9/2$ & $172$ & $D_{7}$ \\ 
     \midrule
    $[4,2,6]$ & $120$ & $11$ & & $5$ & outer & $-2$ & $1$ & \\ 
    &&& & $30$ & outer & $-2$ & $1$ & \\ 
    \cmidrule{4-10}
    &&& $1, 2$ & $6$ & inner & $-3/2$ & $3/2$ & $14$ & $A_{1}+A_{2}+A_{4}$ \\ 
    &&& $3$ & $20$ & inner & $-6/5$ & $3$ & $26$ & $2A_{1}+A_{5}$ \\ 
    &&& $4$ & $30$ & inner & $-6/5$ & $3$ & $26$ & $2A_{1}+A_{5}$ \\ 
    &&& $5$ & $1$ & inner & $-1$ & $2$ & $19$ & $A_{3}+A_{4}$ \\ 
    &&& $6$ & $30$ & inner & $-4/5$ & $4$ & $51$ & $2A_{1}+D_{5}$ \\ 
    &&& $7$ & $40$ & inner & $-4/5$ & $4$ & $51$ & $2A_{1}+D_{5}$ \\ 
    &&& $8$ & $60$ & inner & $-4/5$ & $4$ & $51$ & $2A_{1}+D_{5}$ \\ 
    &&& $9, 10$ & $20$ & inner & $-7/10$ & $7/2$ & $46$ & $A_{1}+A_{6}$ \\ 
    &&& $11, 12$ & $20$ & inner & $-3/10$ & $9/2$ & $146$ & $A_{1}+E_{6}$ \\ 
    &&& $13, 14$ & $60$ & inner & $-3/10$ & $9/2$ & $146$ & $A_{1}+E_{6}$ \\ 
    &&& $15$ & $10$ & inner & $-1/5$ & $4$ & $171$ & $D_{7}$ \\ 
    \midrule
    $[2,0,12 ]$  & $720$ & $15/2$ & & $45$ & outer & $-2$ & $1$ & \\ 
    &&& & $45$ & outer & $-1/2$ & $7/2$ & \\ 
    \cmidrule{4-10}
    &&& $1$ & $10$ & inner & $-3/2$ & $3/2$ & $21/2$ & $2A_{2}+A_{3}$ \\ 
    &&& $2$ & $30$ & inner & $-4/3$ & $2$ & $27/2$ & $A_{1}+2A_{3}$ \\ 
    &&& $3$ & $72$ & inner & $-5/4$ & $5/2$ & $35/2$ & $A_{1}+A_{2}+A_{4}$ \\ 
    &&& $4$ & $60$ & inner & $-1$ & $3$ & $51/2$ & $A_{1}+A_{2}+D_{4}$ \\ 
    &&& $5$ & $12$ & inner & $-5/6$ & $5/2$ & $45/2$ & $A_{3}+A_{4}$ \\ 
    &&& $6$ & $40$ & inner & $-3/4$ & $3$ & $63/2$ & $A_{2}+A_{5}$ \\ 
    &&& $7, 8$ & $120$ & inner & $-7/12$ & $7/2$ & $99/2$ & $A_{1}+A_{6}$ \\ 
    &&& $9$ & $120$ & inner & $-1/3$ & $4$ & $207/2$ & $A_{1}+D_{6}$ \\ 
    &&& $10$ & $180$ & inner & $-1/3$ & $4$ & $207/2$ & $A_{1}+D_{6}$ \\ 
    &&& $11, 12$ & $120$ & inner & $-1/12$ & $4$ & $783/2$ & $E_{7}$ \\ 
    %
    %
    \midrule
    $[6,0,6]$ & $1440$ & $5$ & & $60$ & outer & $-2$ & $1$ & \\ 
    \cmidrule{4-10}
    &&& $1$ & $40$ & inner & $-3/2$ & $3/2$ & $8$ & $A_{1}+3A_{2}$ \\ 
    &&& $2$ & $180$ & inner & $-4/3$ & $2$ & $11$ & $2A_{1}+A_{2}+A_{3}$ \\ 
    &&& $3$ & $10$ & inner & $-1$ & $2$ & $13$ & $2A_{2}+A_{3}$ \\ 
    &&& $4, 5$ & $144$ & inner & $-5/6$ & $5/2$ & $20$ & $A_{1}+A_{2}+A_{4}$ \\ 
    &&& $6$ & $240$ & inner & $-2/3$ & $3$ & $32$ & $2A_{1}+A_{5}$ \\ 
    &&& $7$ & $360$ & inner & $-2/3$ & $3$ & $32$ & $2A_{1}+A_{5}$ \\ 
    &&& $8$ & $180$ & inner & $-1/3$ & $3$ & $59$ & $A_{2}+D_{5}$ \\ 
    &&& $9, 10$ & $240$ & inner & $-1/6$ & $7/2$ & $152$ & $A_{1}+E_{6}$ \\ 
    &&& $11, 12$ & $720$ & inner & $-1/6$ & $7/2$ & $152$ & $A_{1}+E_{6}$  
  \end{longtable}

\end{center}
%
%
\section{Enriques involutions and Borcherds method}\label{sec:EnriquesInvols}
In this section, we assume that $X$ is a complex K3 surface
admitting a primitive embedding $\iota_X\colon S_X\inj L_{26}$ of simple Borcherds type 
and, in addition, that 
\[
\textrm{the natural homomorphism $\rho_X\colon \Aut(X)\to \OG(S_X, \PPP_X)$ 
is injective.} \tag{C} \label{eq:condC}
\]
\subsection{Inner faces}\label{subsec:innerfaces}
Let $D_0$ be an $\Xinduced$-chamber contained in $N_X$.
Let $w_1, \dots, w_k$ be the inner walls of $D_0$.
For each $w_i$,
we calculate an extra automorphism $g_i\in \aut(X)$
associated with $w_i$ (see~Definition~\ref{def:extraX}).
\begin{definition}
A face $f$ of $D_0$ is said to be \emph{$D_0$-inner}
if $f$ is not contained in any outer wall of $D_0$,
whereas
 $f$ is said to be \emph{$N_X$-inner}
if $f$ is not contained in any wall of $N_X$.
\end{definition}
\begin{remark}
An $N_X$-inner face is always $D_0$-inner.
The converse is, however, not true in general
as illustrated in Figure~\ref{fig:example},
in which a black circle indicates a $D_0$-inner face of codimension $2$
that is not $N_X$-inner.
\end{remark}
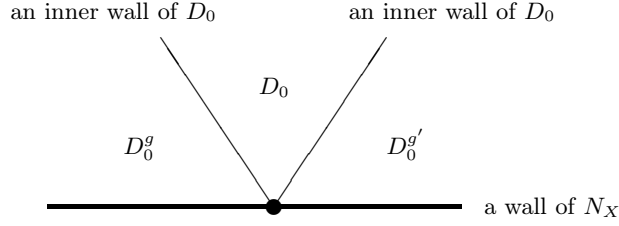
\begin{figure}
{\small
\setlength{\unitlength}{1mm}
\begin{picture}(100, 35)(0, 5)
\put(40, 5){\circle*{2}}
\put(40,5){\line(2, 3){15}}
\put(40,5){\line(-2, 3){15}}
\put(38, 20){$D_0$}
\put(20, 12){$D_0^g$}
\put(55, 12){$D_0^{g\sprime}$}
\put(5, 30){an inner wall of $D_0$}
\put(50, 30){an inner wall of $D_0$}
\linethickness{1.5pt}
\put(10,5){\line(1, 0){55}}
\put(68, 4){a wall of $N_X$}
\end{picture}
}
\caption{A $D_0$-inner face that is not $N_X$-inner.}\label{fig:example}
\end{figure}
%
%
%
%
Let $f$ be a $D_0$-inner face of dimension $>0$.
We put 
\begin{align*}
\DDD(f) &:= \shortset{D}{\textrm{$D$ is an $\Xinduced$-chamber contained in $N_X$ and containing $f$}}, \\
\AAA(X, f) &:=
\shortset{g\in \aut(X)}{D_0^g\in \DDD(f)}=\shortset{g\in \aut(X)}{f\subset D_0^g}, \\
\aut(X, f) &:= \shortset{g\in \aut(X)}{f^g=f}.
\end{align*}
The set $\DDD(f)$ is calculated by the following method.
\begin{algorithm}
We set
$\DDD=[D_0]$, $\gamma_0=\id$, $\varGamma=[\gamma_0]$, and $i=0$.
During the calculation,
the ordered set $\DDD$ is a subset of $\DDD(f)$,
and the $(i+1)$st member 
$\gamma_i$ of $\varGamma$ is an element of $\aut(X)$ that maps $D_0$ to the $(i+1)$st member 
$D_i$ of $\DDD$.
While $i<|\DDD|$,
we execute the following.
We calculate the set 
$\{w_{\nu(1)}, \dots, w_{\nu(m)}\}$ of inner walls $w_{\nu(j)}$ of $D_0$
such that $f\subset w_{\nu(j)}^{\gamma_i}$.
Let $g_{\nu(j)}\in \aut(X)$ be an extra automorphism associated with 
 $w_{\nu(j)}$. %
For each $j=1, \dots, m$, 
we calculate the induced chamber $D\sprime:=D_0^{g_{\nu(j)}\gamma_i}$,
which is adjacent to $D_i=D_0^{\gamma_i}$ across 
$w_{\nu(j)}^{\gamma_i}$ and contains $f$.
If $D\sprime$ has not yet been added to $\DDD$,
we add $D\sprime$ to $\DDD$
and $g_{\nu(j)}\gamma_i$ to $\varGamma$.
Then we increment $i$ to $i+1$.
\end{algorithm}
When this algorithm terminates,
the list $\DDD$ is equal to $\DDD(f)$.
Moreover, we have calculated $\varGamma=\shortset{g_D}{D\in \DDD(f)}$,
where  $g_D\in \aut(X)$ maps $D_0$ to  $D\in \DDD(f)$.
Note that the action of $g_D\in \varGamma$ preserves the walls of $N_X$.
The following is obvious from the definition.
\begin{criterion}
The $D_0$-inner face $f$ is $N_X$-inner
if and only if,
for any $g_D\in \varGamma$
and any outer wall $D_0\cap (r)\sperp$ of $D_0$,
the wall $(D_0\cap (r)\sperp)^{g_D}$ of $D=D_0^{g_D}$
does not contain $f$.
\end{criterion}
Suppose that $f$ is $N_X$-inner and $D$ is an element of $\DDD(f)$.
Note that the set of all elements $g\in \aut(X)$ that maps $D_0$ to $D$
is equal to $\aut(X, D_0)\cdot g_D$.
Therefore we can calculate $\AAA(X, f)$ by 
\[
\AAA(X, f)=\bigsqcup_{D\in \DDD(f)} \aut(X, D_0)\cdot g_D.
\]
The subgroup
$\aut(X, f)$
of $\aut(X)$ is contained in the finite set $\AAA(X, f)$,
and thus we can calculate $\aut(X, f)$.
\begin{definition}
Let $f$ and $f\sprime$ be $N_X$-inner faces of $D_0$.
We say that $f$ and $f\sprime$ are \emph{$\aut(X)$-equivalent}
(resp.~\emph{$\aut(X, D_0)$-equivalent})
if there exists an element $g\in \aut(X)$ 
(resp.~$g\in \aut(X, D_0)$)
such that $f^g=f\sprime$.
\end{definition}
Even though $\aut(X)$ is infinite in general,
we can calculate the $\aut(X)$-equivalence classes by the following:
\begin{criterion}
The faces $f$ and $f\sprime$ are $\aut(X)$-equivalent
if and only if
there exists an element $g\in \AAA(X, f\sprime)$ such that $f^g=f\sprime$.
\end{criterion}
\subsection{An algorithm to classify all Enriques involutions}
Let $\tilde{\enr}\colon X\to X$ be an Enriques involution,
and $\pi\colon X\to Y:=X/\gen{\tilde{\enr}}$
the quotient morphism to the Enriques surface $Y$.
Let $\enr\in \aut(X)$ denote the image of $\tilde{\enr}$ by 
the natural homomorphism~\eqref{eq:nat_rho}.
Then $\pi$ induces a primitive embedding
$\pistar\colon S_Y(2)\inj S_X$.
We have canonical identifications $S_Y(2)\tensor \R=S_Y\tensor \R$
and $\OG(S_Y(2))=\OG(S_Y)$.
In particular,
we regard the positive cone $\PPP_Y$ of $S_Y$ as a positive cone of $S_Y(2)$.
The embedding $\pistar$ induces an embedding
\[
\pistar\colon \PPP_Y\inj \PPP_X.
\]
Henceforth, 
we regard $S_Y(2)$ as a primitive sublattice of $S_X$ and 
$\PPP_Y$ as a subspace of $\PPP_X$ by $\pistar$.
Note that 
$S_Y(2)$ is equal to $\shortset{v\in S_X}{v^\enr=v}$, and 
$\PPP_Y$ is equal to $\shortset{x\in \PPP_X}{x^{\enr}=x}$.
\begin{proposition}\label{prop:NXNY}
We have 
$N_Y=N_X\cap \PPP_Y$.
Let $y$ be a point of $N_Y$. Then 
$y$ is an interior point of $N_Y$
if and only if $y$ is an interior point of $N_X$.
\end{proposition}
\begin{proof}
The first equality is obvious.
Since $\pi$ is \'etale,
the orthogonal complement of $S_Y(2)$ in $S_X$ contains no $(-2)$-vectors,
and a line bundle of $Y$ is ample if and only if
its pull-back to $X$ is ample.
\end{proof}
Let $y$ be a sufficiently general point of $N_Y$.
By Theorem~\ref{thm:D0}, there exists an automorphism $g\in \aut(X)$ such that
$y^g\in D_0$, and hence 
$D_0\cap N_Y^g $ contains a non-empty open subset of $\PPP_Y^g$.
Therefore, replacing $\enr$ by $g\inv \enr g$,
we can assume that
\[
E_0:= D_0\cap N_Y
\]
contains a non-empty open subset of $\PPP_Y$.
Consider the composite
\[
\iotaY:=\iotaX\circ \pistar \colon S_Y(2)\inj L_{26}
\]
of primitive embeddings.
Then $\PPP_Y$ is decomposed into the union of $\iotaY^*\RRR_{26}\sperp$-chambers.
Since every wall of $N_Y$ is defined by a $(-2)$-vector, 
it follows that $N_Y$ is decomposed into a union of
$\Yinduced$-chambers.
Note that 
$E_0$ is one of the $\Yinduced$-chambers in $N_Y$.
\begin{definition}
For a closed subset $A$ of $D_0$,
the \emph{minimal face} of $D_0$ for $A$ is 
the face of $D_0$ containing $A$ with the minimal dimension.
\end{definition}

Let $f_{\enr}$ be the minimal face of $D_0$ for $E_0$.
Since 
the orthogonal complement of $S_Y(2)$ in $S_X$ contains no $(-2)$-vector,
the face $f_{\enr}$ is $N_X$-inner.
Moreover, the involution $\enr\in \aut(X)$ belongs to
$\aut(X, f_{\enr})$.
Let $\enr\sprime$ be 
an Enriques involution 
such that $f_{\enr\sprime}$ is a face of $D_0$.
If $\enr\sprime$ 
is conjugate to $\enr$, 
then $f_{\enr}$ is $\aut(X)$-equivalent to $f_{\enr\sprime}$.
If $f_{\enr}=f_{\enr\sprime}$,
then $\enr$ and $\enr\sprime$ are conjugate 
if and only if $\enr$ and $\enr\sprime$
are conjugate in $\aut(X, f_{\enr})$.

We calculate all $N_X$-inner faces of $D_0$ of dimension $\ge 10$ 
by descending induction of 
the dimension of faces (see~Section~\ref{subsec:defchambersfaces}),
and compute a complete set of representatives of 
the $\aut(X)$-equivalence classes.
For each representative $f$,
we calculate $\aut(X, f)$.
We then calculate the set of Enriques involutions $\enr$ contained in $\aut(X, f)$
such that $f_{\enr}=f$
by Keum's criterion (Theorem~\ref{thm:keum}),
and thus 
we obtain a set of complete representatives of
Enriques involutions in $\aut(X)$ modulo conjugation.
\subsection{Computation of \texorpdfstring{$\Aut(Y)$}{Aut(Y)}}
Let $\enr$ be a representative
of $\aut(X)$-conjugacy classes of Enriques involutions 
obtained by the method above.
In particular, we have an $\Yinduced$-chamber
$E_0=D_0\cap N_Y$, the minimal face $f_{\enr}$ of $D_0$ for $E_0$,
and the associated data $\DDD(f_{\enr})$, $\AAA(X, f_{\enr})$, $\aut(X, f_{\enr})$.
We put
\[
\aut(X, \enr):=\shortset{g_X\in \aut(X)}{\enr g_X=g_X\enr}=\shortset{g_X\in \aut(X)}{S_Y(2)^{g_X}=S_Y(2)},
\]
where the second equality follows from $S_Y(2)=\shortset{v\in S_X}{v^\enr=v}$.
We have a natural restriction homomorphism
$\aut(X, \enr)\to \OG(S_Y)$,
which is denoted by $g_X\mapsto g_X|S_Y$.
By Condition \eqref{eq:condC},
we have a natural identification
\begin{equation}\label{eq:AutYcen}
\Aut(Y)\cong\aut(X, \enr)/\gen{\enr}.
\end{equation}
Under the identification~\eqref{eq:AutYcen},
the homomorphism $\rho_Y\colon \Aut(Y)\to \OG(S_Y, \PPP_Y)$ 
is identified with the  homomorphism
$g_X\bmod \gen{\enr}\mapsto g_X|S_Y$.
The method below,
when it works,
gives us a finite set of generators of $\aut(X, \enr)$,
and hence a finite set of generators of $\Aut(Y)$.
\par
Recall that $\aut(Y)$ is the image of $\Aut(Y)$ by $\rho_Y$. 
We put
\[
\aut(Y, E_0):=\shortset{g\in \aut(Y)}{E_0^g=E_0},
\]
and let $\Aut(Y, E_0)$ denote the inverse image of $\aut(Y, E_0)$ by $\rho_Y$.

\begin{proposition}\label{prop:autYpreservesinduced}
The action of $\aut(Y)$ on $N_Y$
preserves the set of $\Yinduced$-chambers contained in $N_Y$.
\end{proposition}
\begin{proof}
Let $g$ be an element of $\aut(Y)$.
Then $g$ extends to $g_X\in \aut(X, \enr)$.
By Condition \eqref{eq:condA},
this isometry $g_X\in \OG(S_X, \omega_X) \cap \OG(S_X,\cP_X)$
extends to an isometry $\tilde{g}_X$ of $L_{26}$,
which preserves the set of Conway chambers.
Hence its restriction $g$ to $S_Y(2)$
preserves the set of chambers induced by Conway chambers.
\end{proof}
We put
\[
\aut(X, \enr, f_{\enr}):=\aut(X, \enr)\cap \aut(X, f_{\enr}).
\]
\begin{proposition}\label{prop:AutYE0}
The identification~\eqref{eq:AutYcen} induces
$\Aut(Y, E_0)\cong \aut(X, \enr, f_{\enr})/\gen{\enr}$.
\end{proposition}
\begin{proof}
Note that $E_0=f_{\enr}\cap N_Y$.
Since $E_0$ contains an interior point of 
the face $f_{\enr}$,
an element $g_X$ of $\aut(X, \enr)$
fixes $E_0$ if and only if $g_X$ fixes $f_{\enr}$.
\end{proof}
\begin{corollary}\label{cor:KerrhoY}
By the identification~\eqref{eq:AutYcen},
the kernel of $\rho_Y\colon \Aut(Y)\to \OG(S_Y, \PPP_Y)$ is equal to 
\[
\shortset{g_X\in \aut(X, \enr, f_{\enr})}{g_X|S_Y=\id}/\gen{\enr}.
\]
\end{corollary}
Recall from Section~\ref{subsec:LtttoLts} that we have classified 
primitive embeddings of $S_Y(2)\cong L_{10}(2)$ into $L_{26}$.
%
%
The  $\Yinduced$-chamber $E_0$ has only finitely many walls.
By~Remark~\ref{rem:infty},
the primitive embedding $\iota_Y\colon S_Y(2)\inj L_{26}$
is not of type \BStype{infty}.
By Theorem~\ref{thm:BS}, 
every $\Yinduced$-chamber $E$ has only a finite number of walls,
and each wall of $E$ is defined by a $(-2)$-vector $r\in \RRR_Y$.
\begin{definition}
A wall $w$ of $E_0$ is said to be \emph{outer} 
if $w$ is contained in a wall of~$N_Y$.
Otherwise $w$ is said to be \emph{inner}.
\end{definition}
There are several criteria 
to determine whether a given wall $w$ of $E_0$ is outer or inner.
\begin{criterion}
Suppose that the wall $w$ of $E_0$ is defined by $r\in \RRR_Y$.
Then $w$ is outer if and only if there exists a $(-2)$-vector 
$u$ in the orthogonal complement $[\pi^*]\sperp$ of $S_Y(2)$ in $S_X$
such that $(u+r)/2\in S_X$.
\end{criterion}
Indeed, 
the condition in the statement is equivalent to the condition 
that $r$ is the class of an effective divisor of $Y$
(see~\cite{Nikulin1984}).
\begin{criterion}
Let $f_{\enr}(w)$ be the minimal face of $D_0$ for the closed subset $w$ of~$D_0$.
Then $w$ is inner if and only if $f_{\enr}(w)$ is $N_X$-inner.
\end{criterion}
Indeed, 
by minimality of $f_{\enr}(w)$,
there exists an interior point $y$ of $w$ that is an interior point of $f_{\enr}(w)$.
Then the statement follows from Proposition~\ref{prop:NXNY}.
%
%
%
\par
When $E_0$ has no inner walls,
we have $E_0=N_Y$ and 
$|\Aut(Y)|<\infty$, and 
the Nikulin-Kondo type 
of $Y$
is obtained by comparing the configuration of
$(-2)$-vectors defining the walls of $E_0$
with the dual graphs of smooth rational curves given in~\cite{Kondo1986}.
\par
We consider $\Aut(Y)$ when $E_0$ has an inner wall.
Let $I_0$ denote the set of inner walls of $E_0$.
For each $w=E_0\cap (r)\sperp\in I_0$ with $r\in \RRR_Y$,
we put $E(w):=E_0^{s_r}$,
where $s_r\colon \PPP_Y\to \PPP_Y$ is the reflection into the hyperplane $(r)\sperp\subset \PPP_Y$.
Theorem~\ref{thm:BS} implies that 
$E(w)$ is the $\Yinduced$-chamber adjacent to $E_0$ across $w$.
Recall that $\AAA(X, f_{\enr}(w))$ is the set of $g_X\in \aut(X)$
such that $D_0^{g_X}$ contains $f_{\enr}(w)$.
If the restriction $g_X|S_Y$ to $S_Y(2)$ 
of $g_X\in \aut(X, \enr)$ maps $E_0$ to $E(w)$, then 
$g_X\in \AAA(X, f_{\enr}(w))$ holds.
\begin{definition}
An element $g_X$ of $\aut(X, \enr)\cap \AAA(X, f_{\enr}(w))$
is an \emph{extra automorphism} for the inner wall $w\in I_0$
if the restriction  $g_X|S_Y$ of $g_X$ to $S_Y(2)$ maps $E_0$ to 
 $E(w)$.
\end{definition}
Since $\AAA(X, f_{\enr}(w))$ is finite,
we can determine the existence of an extra automorphism
for each inner wall of $E_0$.
\begin{theorem}\label{thm:cen}
Suppose that Condition~\eqref{eq:condC} is satisfied.
Suppose also that the following holds:
\[
 \textrm{there exists an extra automorphism $g_X(w)$ for each inner wall $w\in I_0$.} \tag{IY} \label{eq:condIY}
\]
Then $\aut(X, \enr)$
is generated by the finite subgroup $\aut(X, \enr, f_{\enr})$
and the extra automorphisms $g_X(w)$ {\rm ($w\in I_0$)}.
\end{theorem}
\begin{proof}
Let $\Gamma$ denote the subgroup of
$\aut(X, \enr)$ generated by the extra automorphisms $g_X(w)$ {\rm ($w\in I_0$)}.
First we prove the following claim.
For any $\Yinduced$-chamber~$E$ contained in $N_Y$,
there exists an element $\gamma\in \Gamma$
such that $\gamma|S_Y$ maps $E_0$ to $E$.
There exists a chain 
$E_0, E_1, \dots, E_m=E$
of $\Yinduced$-chambers contained in $N_Y$
such that $E_{i-1}$ and $E_{i}$ is adjacent
for $i=1, \dots, m$.
We prove the claim by induction on 
the length $m$ of the chain
with the case $m=0$ being trivial.
Suppose that $m>0$.
There exists an element $\gamma\sprime\in \Gamma$
such that $\gamma\sprime|S_Y$ maps $E_0$ to $E_{m-1}$.
Let $E\sprime$ be the $\Yinduced$-chamber
that is mapped to $E_m$ by $\gamma\sprime|S_Y$.
Then $E\sprime$ is  adjacent to $E_0$.
Note that $\gamma\sprime|S_Y\in \aut(Y)$ preserves $N_Y$.
Therefore $E\sprime$ is contained in $N_Y$.
In particular, 
the wall $w$ between $E_0$ and $E\sprime$ is inner,
and hence there exists an extra automorphism $g_X(w)$
such that $g_X(w)|S_Y$ maps $E_0$ to $E\sprime$.
We put $\gamma:=g_X(w)\cdot \gamma\sprime\in \Gamma$.
Then $\gamma|S_Y$ maps $E_0$ to $E_m$.
\par
Next we show that $\Gamma$ and $\aut(X, \enr, f_{\enr})$
generate $\aut(X, \enr)$.
Let $g$ be an arbitrary element of $\aut(X, \enr)$.
We apply the claim above to the $\Yinduced$-chamber $E_0^{g|S_Y}$,
and obtain an element $\gamma\in \Gamma$
such that $(g\gamma\inv)|S_Y$ is an element of $\aut(Y, E_0)$.
By Proposition~\ref{prop:AutYE0}, we have $g\gamma\inv\in \aut(X, \enr, f_{\enr})$.
\end{proof}
%
%
\begin{definition}\label{def:simpleenr}
We say that a triple $(X, \iota_X, \enr)$
of a K3 surface $X$, a primitive embedding $\iota_X\colon S_X\inj L_{26}$,
and an Enriques involution $\enr$ of $X$ 
is of \emph{simple Borcherds type}
if $X$ satisfies Condition~\eqref{eq:condC}, 
$(X, \iota_X)$ is of simple Borcherds type
in the sense of Definition~\ref{def:simple}, 
and $\enr$ satisfies Condition~\eqref{eq:condIY}.
\end{definition}
\begin{remark}
The notion of simple Borcherds type was introduced in~\cite{ShimadaHoles}
for K3 surfaces.
We hope that we can find 
a bound on the degrees of polarizations similar to that of~\cite{ShimadaHoles}
for Enriques surfaces.
\end{remark}
\subsection{Enriques involutions of the 11 singular K3 surfaces} \label{sec:enriques-involutions-of-certain-K3}
\begin{table}
\small
\caption{Enriques involutions of the 11 singular K3 surfaces.}
\label{table:enrtable}
\begin{tabular}{ccccccccccc}
\toprule
 No. & $T_X$ & $\dim f_{\enr}$ & $\iota_Y$ & ${\rm NK}$ & \texttt{m4} & $|{\rm ws}|$ & $|G_{\enr}|$ & $|I_0|$ & $|K_{\rho}|$ & $|\mathrm{aut}|$ \\ 
\midrule 
1 & $[2, 1, 4]$ & $19$ & $\BStype{12B}$ & $\NKtype{2}$ & $144$ & $12$ & $48$ & $0$ & $1$ & $24$ \\ 
2 &  & $18$ & $\BStype{12A}$ & $\NKtype{1}$ & $242$ & $12$ & $16$ & $0$ & $2$ & $4$ \\ 
\midrule 
3 & $[2, 0, 6]$ & $19$ & $\BStype{12B}$ & $\NKtype{2}$ & $144$ & $12$ & $48$ & $0$ & $1$ & $24$ \\ 
\midrule 
4 & $[4, 2, 4]$ & $18$ & $\BStype{12A}$ & $\NKtype{1}$ & $246$ & $12$ & $16$ & $0$ & $2$ & $4$ \\ 
5 &  & $18$ & $\BStype{20B}$ & $\NKtype{3}$ & $246$ & $20$ & $64$ & $4$ & $2$ & $\infty$ \\ 
6 &  & $17$ & $\BStype{20A}$ & $\NKtype{5}$ & $246$ & $20$ & $96$ & $0$ & $2$ & $24$ \\ 
\midrule 
7 & $[2, 1, 8]$ & $19$ & $\BStype{20D}$ & $\NKtype{7}$ & $90$ & $20$ & $120$ & $5$ & $1$ & $\infty$ \\ 
8 &  & $19$ & $\BStype{12B}$ & $\NKtype{2}$ & $144$ & $12$ & $48$ & $0$ & $1$ & $24$ \\ 
9 &  & $19$ & $\BStype{12B}$ & $\NKtype{2}$ & $144$ & $12$ & $48$ & $0$ & $1$ & $24$ \\ 
10 &  & $18$ & $\BStype{12A}$ & $\NKtype{1}$ & $240$ & $12$ & $8$ & $2$ & $2$ & $\infty$ \\ 
11 &  & $17$ & $\BStype{20A}$ & $\NKtype{5}$ & $132$ & $20$ & $48$ & $4$ & $1$ & $\infty$ \\ 
\midrule 
12 & $[4, 0, 4]$ & $20$ & $\BStype{20F}$ & $\NKtype{4}$ & $180$ & $20$ & $640$ & $0$ & $1$ & $320$ \\ 
13 &  & $19$ & $\BStype{20D}$ & $\NKtype{7}$ & $180$ & $20$ & $120$ & $5$ & $1$ & $\infty$ \\ 
14 &  & $19$ & $\BStype{12B}$ & $\NKtype{2}$ & $180$ & $12$ & $48$ & $0$ & $1$ & $24$ \\ 
15 &  & $18$ & $\BStype{12A}$ & $\NKtype{1}$ & $244$ & $12$ & $16$ & $0$ & $2$ & $4$ \\ 
16 &  & $18$ & $\BStype{12A}$ & $\NKtype{1}$ & $244$ & $12$ & $16$ & $2$ & $4$ & $\infty$ \\ 
17 &  & $18$ & $\BStype{20B}$ & $\NKtype{3}$ & $244$ & $20$ & $64$ & $8$ & $2$ & $\infty$ \\ 
18 &  & $18$ & $\BStype{20B}$ & $\NKtype{3}$ & $244$ & $20$ & $64$ & $4$ & $2$ & $\infty$ \\ 
19 &  & $18$ & $\BStype{20B}$ & $\NKtype{3}$ & $308$ & $20$ & $256$ & $0$ & $2$ & $64$ \\ 
20 &  & $17$ & $\BStype{20A}$ & $\NKtype{5}$ & $244$ & $20$ & $32$ & $4$ & $2$ & $\infty$ \\ 
\midrule 
21 & $[4, 2, 6]$ & $19$ & $\BStype{20D}$ & $\NKtype{7}$ & $92$ & $20$ & $240$ & $0$ & $1$ & $120$ \\ 
22 &  & $18$ & $\BStype{12A}$ & $\NKtype{1}$ & $242$ & $12$ & $16$ & $0$ & $2$ & $4$ \\ 
\midrule 
23 & $[2, 0, 12]$ & $19$ & $\BStype{20D}$ & $\NKtype{7}$ & $90$ & $20$ & $120$ & $5$ & $1$ & $\infty$ \\ 
\midrule 
24 & $[6, 0, 6]$ & $20$ & $\BStype{40E}$ & $ $ & $60$ & $40$ & $1440$ & $10$ & $1$ & $\infty$ \\ 
25 &  & $18$ & $\BStype{12A}$ & $\NKtype{1}$ & $240$ & $12$ & $16$ & $2$ & $4$ & $\infty$ \\ 
26 &  & $17$ & $\BStype{20A}$ & $\NKtype{5}$ & $132$ & $20$ & $48$ & $4$ & $1$ & $\infty$ \\ 
\bottomrule
\end{tabular}
\end{table}
We apply the method in the previous section 
to the singular K3 surfaces in Section~\ref{subsec:singK3}.
First remark that 
Condition (C) holds for the $11$ cases except for the cases $T_X=[2,1,2]$ and $T_X=[2,0,2]$
(see Remark~\ref{rem:Kerrho} and Table~\ref{table:Table1}).
Note that in these two cases, and also in the case $T_X = [2,0,4]$, there exist no Enriques involutions by Theorem~\ref{thm:sertoez}.
\par
Our main result is as follows.
%
\begin{theorem}\label{thm:ichiro-main}
Let $X$ be one of the singular K3 surfaces of No.~$\ne 1, 2, 4$
in Table~\ref{table:Table1},
and let $\iota_X\colon S_X\inj L_{26}$ be the primitive embedding given in 
Section~\ref{subsec:singK3}.
Then the Enriques involutions of $X$ modulo conjugation in $\Aut(X)\cong \aut(X)$
are given in Table~\ref{table:enrtable}.
For each Enriques involution $\enr$ on $X$,
the triple $(X, \iota_X, \enr)$ is of simple Borcherds type.
\end{theorem}
We explain the contents of Table~\ref{table:enrtable}.
The item $\iota_Y$ is the type of the primitive embedding $\iota_Y\colon S_Y(2)\inj L_{26}$
given in~\cite{BS}.
The item {\rm NK} is
the Nikulin-Kondo type
of the $\Yinduced$-chamber $E_0$ (see Theorem~\ref{thm:sigmatau}).
The item \texttt{m4} is the number of $(-4)$-vectors in 
the orthogonal complement of $S_Y(2)$ in $S_X$.
The item $|{\rm ws}|$
is the number of walls of $E_0$.
The item $|G_{\enr}|$
is the order of 
\[
G_{\enr}:=\aut(X, \enr, f_{\enr}).
\]
The item $|I_0|$ is the number of inner walls of $E_0$.
\begin{remark}
For the Enriques involution No.\,24 on $X$
with $T_X=[6,0,6]$,
the $\Yinduced$-chamber $E_0$ has $40$ walls
and the configuration of the walls is not of Nikulin-Kondo type.
The dual graph 
is too complicated to be presented here.
See~\cite{thecompdata} for the matrix presentation of this configuration.
\end{remark}
The item $|K_{\rho}|$ is the order of 
the kernel of $\rho_Y\colon \Aut(Y)\to \aut(Y)$,
and the item $|\mathrm{aut}|$ is the order of $\aut(Y)$.
The fact that $\aut(Y)$ is infinite when $I_0$ is non-empty
was confirmed by selecting elements of $\aut(Y)$
randomly by means of the finite generating set of $\aut(Y)$ obtained by~Theorem~\ref{thm:cen}
and finding a matrix of infinite order among these sample elements.
\begin{figure}
\begin{center}
{\footnotesize
\setlength{\unitlength}{.64mm}
\begin{picture}(120, 40)
\put(10, 20){\circle{6}}
\put(30, 20){\circle{6}}
\put(50, 20){\circle{6}}
\put(70, 20){\circle{6}}
\put(90, 20){\circle{6}}
\put(110, 20){\circle{6}}
\put(10, 5){\circle{6}}
\put(60, 5){\circle{6}}
\put(110, 5){\circle{6}}
\put(10, 35){\circle{6}}
\put(60, 35){\circle{6}}
\put(110, 35){\circle{6}}
\put(8.8, 18.5){$3$}
\put(29, 18.5){$9$}
\put(47.7, 18.35){$11$}
\put(67.7, 18.35){$12$}
\put(87.7, 18.35){$10$}
\put(109, 18.5){$7$}
\put(8.8, 3.5){$4$}
\put(59, 3.5){$5$}
\put(109, 3.5){$6$}
\put(8.8, 33.5){$2$}
\put(59, 33.5){$1$}
\put(109, 33.5){$8$}
\put(13,20){\line(1, 0){14}}
\put(33,21){\line(1, 0){14}}
\put(33,19){\line(1, 0){14}}
\put(53,21){\line(1, 0){14}}
\put(53,19){\line(1, 0){14}}
\put(73,21){\line(1, 0){14}}
\put(73,19){\line(1, 0){14}}
\put(93,20){\line(1, 0){14}}
\put(13,5){\line(1, 0){44}}
\put(63,5){\line(1, 0){44}}
\put(13,35){\line(1, 0){44}}
\put(63,35){\line(1, 0){44}}
\put(10,23){\line(0, 1){9}}
\put(10,17){\line(0, -1){9}}
\put(110,23){\line(0, 1){9}}
\put(110,17){\line(0, -1){9}}
\end{picture}
}
\end{center}
\caption{Configuration of Nikulin-Kondo type $\NKtype{1}$}\label{fig:I}
\end{figure}
\begin{remark}
Consider the Enriques involutions of Nos.~10,~16 and~25,
that is, the cases where the Nikulin-Kondo type is $\NKtype{1}$ and $\Aut(Y)$ is infinite.
In these cases, we have $|I_0|=2$.
The configuration of Nikulin-Kondo type $\NKtype{1}$ is as in Figure~\ref{fig:I},
and the inner walls are defined by the $(-2)$-vectors \textcircled{\scriptsize 11} and \textcircled{\scriptsize 12}.
\par
See~\cite{thecompdata}
for the inner walls of $E_0$ for the other Enriques involutions. 
The finite generating sets of $\aut(X, \enr)$ 
and of $\aut(Y)$ are also given explicitly in~\cite{thecompdata}.
\end{remark}

Table~\ref{tab:Facedata_all} is a list of $N_X$-inner faces of $D_0$
that corresponds to Enriques involutions. 
Note that an $\aut(X)$-equivalence class of $N_X$-inner faces 
is a union 
of orbits of the action of $\aut(X, D_0)$
on the set of $N_X$-inner faces.


\begin{table}
\small

  \caption{$N_X$-inner faces  corresponding to Enriques involutions.}
  \label{tab:Facedata_all}
  \begin{tabular}{ccllccc}
    \toprule
    $T_X$ &  $\dim$ & \texttt{numb} &  \texttt{pws} & $|\DDD|$ & $|\aut(X, f)|$ & $\enr$ \\
    \midrule



    
    $[2,1,4]$ &     $19$ & $14$ & $1^1$ & $2$ & $48$ & No.\,$1$ \\ 
   & $18$ & $42+84$ & $1^2,1^12^1$ & $6$ & $16$ &  No.\,$2$ \\ 
    \midrule
    $[2,0,6]$ &  $19$& $6$ & $3^1$ & $2$ & $48$ & No.\,$3$\\ 
    \midrule
    $[4,2,4]$ & $18$ & $288\times2$ & $1^12^1,1^13^1$ & $8$ & $16$ & No.\,$4$\\ 
    & $18$ & $12$ & $1^2$ & $6$ & $576$ & No.\,$5$\\
    & $17$ & $144$ & $1^22^1$ & $12$ & $96$ & No.\,$6$\\
    \midrule
    $[2,1,8]$ & $19$ & $12$ & $1^1$ & $2$ & $120$ & No.\,$7$\\ 
    & $19$ & $30$ & $4^1$ & $2$ & $48$ & No.\,$8$\\ 
    & $19$ & $30$ & $5^1$ & $2$ & $48$ & No.\,$9$\\ 
    & $18$  & $180\times4$ & $1^14^1,1^15^1,3^14^1,3^15^1$ & $8$ & $8$ & No.\,$10$\\ 
    & $17$  & $90\times2$ & $1^23^1$ & $12$ & $48$ & No.\,$11$\\ 
    \midrule
    $[4,0,4]$ & $20$ & $1$ &  & $1$ & $3840$ & No.\,$12$\\ 
    & $19$ & $64$ & $1^1$ & $2$ & $120$ & No.\,$13$\\ 
    & $19$ & $160$ & $3^1$ & $2$ & $48$ & No.\,$14$\\ 
    & $18$ & $960\times2$ & $1^12^1,1^14^1$ & $8$ & $16$ & No.\,$15$\\ 
    & $18$ & $960\times2$ & $2^13^1$ & $8$ & $16$ & No.\,$16$\\ 
    & $18$ & $60$ & $2^2$ & $4$ & $256$ & Nos.\,$17,18,19$\\ 
    & $17$ & $480+960$ & $1^22^1,1^12^2$ & $12$ & $32$ & No.\,$20$\\ 
    \midrule
    $[4,2,6]$ & $19$ & $1$ & $5^1$ & $2$ & $240$ & No.\,$21$ \\ 
    & $18$ & $30\times2$ & $4^15^1,4^115^1$ & $8$ & $16$ & No.\,$22$ \\ 
    \midrule
    $[2,0,12]$ & $19$ & $12$ & $5^1$ & $2$ & $120$ & No.\,$23$ \\
    \midrule
    $[6,0,6]$ &  $20$ & $1$ &  & $1$ & $1440$ & No.\,$24$ \\ 
     &  $18$ & $360\times2$ & $7^18^1$ & $8$ & $16$ & No.\,$25$ \\
     &  $17$ & $180\times2$ & $2^28^1$ & $12$ & $48$ & No.\,$26$ \\ 
    \bottomrule
    \end{tabular}
\end{table}


The item \texttt{numb}
gives the number of faces in the $\aut(X)$-equivalence class.
The formula in this column shows the decomposition of 
the $\aut(X)$-equivalence class
into a union of $\aut(X, D_0)$-orbits.
The item \texttt{pws} indicates the types of  inner walls of $D_0$ passing
through the face.
The type of an inner wall of $D_0$
is given by No.~in 
 Table~\ref{tab:D0data}.

For example, take the case $T_X = [2,1,4]$.
For a face $f$ in the $\aut(X)$-equivalence class corresponding to the Enriques involution No.~2,
there exist exactly two inner walls of $D_0$
passing through $f$,
and they are both of type 1,
whereas for 
another face $f\sprime$ in this $\aut(X)$-equivalence class,
there exist exactly two inner walls of $D_0$
passing through $f\sprime$,
and they are of type 1 and 2.

\begin{figure}
\begin{center}
{\small
\setlength{\unitlength}{.8mm}
\begin{picture}(100, 58)(0, -20)
\put(40, 5){\circle*{2}}
\put(40,5){\line(2, 3){15}}
\put(40,5){\line(-2, 3){15}}
\put(40,5){\line(-2, -3){15}}
\put(40,5){\line(2, -3){15}}
\put(38, 20){$D\spar{2}$}
\put(20, 12){$D\spar{3}$}
\put(55, 12){$D\spar{1}$}
\put(38, -14){$D\spar{-2}$}
\put(20, -6){$D\spar{-1}$}
\put(55, -6){$D\spar{-3}$}
\put(20, 30){$(v_2)\sperp$}
\put(52, 30){$(v_1\sprime)\sperp$}
\put(10,5){\line(1, 0){55}}
\put(68, 4){$(v_1)\sperp$}
\end{picture}
}
\end{center}
\caption{The $N_X$-inner face $f$.}\label{fig:f5}
\end{figure}
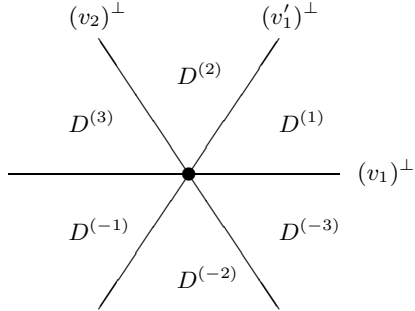

We explain how 
the data \texttt{pws} depends on the choice of a representative of 
an $\aut(X)$-equivalence class.
Let $f$ be a face in this $\aut(X)$-equivalence class.
Then there exist exactly three members 
$(v_1)\sperp, (v_1\sprime)\sperp, (v_2)\sperp$ in 
the family $\Xinduced$ of hyperplanes 
that pass through $f$,
where $v_1$, $v_1\sprime, v_2$ are primitive vectors of $S_X\dual$
such that 
$\intf{v_1, v_1}=\intf{v_1\sprime, v_1\sprime}=-8/7$ and $\intf{v_2, v_2}=-4/7$.
See Figure~\ref{fig:f5}.
If $D_0$ is located in the region $D\spar{\pm 1}$,
then the data \texttt{pws} for $f$ is $1^2$,
whereas if $D_0$ is located in the region $D\spar{\pm 2}$ or $D\spar{\pm 3}$,
then the data \texttt{pws} for $f$ is $1^1 2^1$.
\par
The item $|\DDD|$ is the size of $\DDD(f)$
and $|\aut(X, f)|$ is the order of the group $\aut(X, f)$.
The item $\enr$ shows the Nos. of the  Enriques involutions given in Table~\ref{table:enrtable}.

\section{The two most algebraic Enriques surfaces} \label{sec:themostalgebraic}

In this section,
we study the two most algebraic Enriques surfaces,
that is, Enriques surfaces covered by the singular K3 surface $X_7$ of discriminant $7$.

\subsection{Conjugagy classes of Enriques involutions}

We exemplify Theorem~\ref{thm:davide-main} for the case $X_7$.
Let $T = T_{X_7} = [2,1,4]$ and $S = S_{X_7}$. Let $\iota \in I_{X_7}$ and put $N := [\iota]^\perp(-1)$. 
Let $q := q(T) = \anglefrac{2}{7}$,
so that $q(S) \cong -q \cong \anglefrac{6}{7}$. 
In the notation of Proposition~\ref{prop:Nikulin1.15}, 
the subgroup $H \subset q([\iota])$ must be trivial, 
so $N$ is an even lattice of genus $\fg(10,0,u_{1}^{\oplus 5} \oplus q)$.
By Lemma~\ref{lem:divisible-lattice}, $N \cong N'(2)$, with $N'$ an even lattice of genus $\fg(10,0,q)$.

\begin{lemma} \label{lem:N10_7}
The genus $\fg(10,0,q)$ contains exactly two isomorphism classes, namely $N^{242}_{10,7}$ and $N^{144}_{10,7}$ (see Table~\ref{tab:lattices}).
\end{lemma}
\proof
Let $N'$ be a lattice in this genus. The smallest lattice with bilinear form $b = -b(q)$ is the odd lattice $M_{3,7} := [2,1,2,1,1,3]$, which is unique in its genus. Thus, by \cite{Nikulin1979}, $N' \cong [\iota]^\perp$ for some primitive embedding $\iota\colon M_{3,7} \hookrightarrow L$ into a unimodular lattice $L$ of rank~$13$. Inspecting all such embeddings, we find exactly two non-isomorphic even orthogonal complements, namely 
$N^{242}_{10,7}$ and $N^{144}_{10,7}$.
\endproof

By Proposition~\ref{prop:Nikulin1.5.1}, for both $N = N^{242}_{10,7}(2)$ and $N = N^{144}_{10,7}(2)$, the set $I_{X_7}(N)$ has exactly one $\OG(S)$-orbit. Thus, $r = 2$ in Theorem~\ref{thm:davide-main}. 
Since $\OG(q(S),\omega_{X_7}) = \OG(q(S))$, there is exactly $1$ double coset in both cases. Hence, $X_7$ admits exactly two Enriques involutions up to conjugation in $\aut(X)$.
The two involutions can be distinguished by the number of $(-4)$-vectors in the orthogonal complements of their fixed lattices.

\subsection{Models of the two Enriques quotients}

By the results of Section~\ref{sec:enriques-involutions-of-certain-K3}, the two quotients $Y_\NKtype{1}$ and $Y_\NKtype{2}$ of $X_{7}$ have Nikulin-Kondo type \NKtype{1} and \NKtype{2}. 
Kondo~\cite{Kondo1986} gives two explicit $1$-dimensional families containing all Enriques surfaces of Nikulin-Kondo type \NKtype{1} and \NKtype{2}. Each family depend on one parameter $\alpha$; in this section we determine which values of $\alpha$ give $Y_\NKtype{1}$ and $Y_\NKtype{2}$.
We first summarize Kondo's construction, which is originally due to Horikawa~(\cite{horikawaI}, \cite{horikawaII}; see also Section V.23 in~\cite{bhpv}).

Let $\phi$ be the involution on $\IP^1\times \IP^1$ defined by
\[
 ([u_0,u_1],[v_0,v_1]) \mapsto ([u_0,-u_1],[v_0,-v_1]),
\]
and consider the curves
\[
 L_{1}\colon u_{0} = u_{1}; \quad L_2\colon u_{0} = -u_{1}; \quad L_3\colon v_{0} = v_{1}; \quad L_4\colon v_{0} = -v_{1}.
\]
Let $C$ be a curve of bidegree $(2,2)$, defined by a polynomial $f(u_{0},u_{1},v_{0},v_{1})$, which is invariant with respect to $\phi$, and consider the divisor $B = C + \sum_{i=1}^4 L_i$. 

Let $X$ be the minimal resolution of the double covering of $\IP^1\times \IP^1$ ramified over~$B$. In Kondo's families, $C$ is chosen so that $X$ is a K3 surface and $\phi$ lifts to an Enriques involution $\tilde\enr$ of $X$. We let $Y$ be the quotient of $X$ by $\tilde\enr$.

We wish to find a model of the Enriques quotient $Y$ as an \emph{Enriques sextic surface}, i.e. a non-normal surface of degree $6$ in $\IP^3$ that passes doubly through the edges of the coordinate tetrahedron. 
Such a model is given by a linear system of the form $|F_1 + F_2 + F_3|$, where $F_i$ are half-pencils on $Y$ such that $\langle F_i,F_j \rangle = 1$ (see~\cite{dolgachev-brief-introduction}).

For $i = 1,2$, the composite morphism
\[
 \pi_i\colon X \rightarrow \IP^1 \times \IP^1 \xrightarrow{\pr_i} \IP^1
\]
is an elliptic fibration on $X$, which induces an elliptic fibration $\bar \pi_i$ on $Y$ such that the following diagram commutes
\[
 \begin{tikzcd}
 X \arrow[r,"\pi_i"] \arrow[d] & \IP^1 \arrow[d, "\sigma"] \\
 Y \arrow[r,"\bar\pi_i"] & \IP^1 
 \end{tikzcd}
\]
where $\sigma$ is the map $[u_0,u_1] \mapsto [u_0^2,u_1^2]$ if $i = 1$, or $[v_0,v_1] \mapsto [v_0^2,v_1^2]$ if $i = 2$. The half-pencils on $Y$ are the inverse images of $[0,1],[1,0] \in \IP^1$ (see Section~VIII.17 in~\cite{bhpv}).

There is a third elliptic fibration $\pi_3\colon X \rightarrow \IP^1$, one of whose fiber is the strict transform of~$C$ on~$X$. The half pencils of $\bar\pi_3\colon Y \rightarrow \IP^1$ correspond to the fibers over~$C$ and over $\sum L_i$. We choose coordinates $w_0,w_1$ on $\IP^1$ so that the half-pencils are mapped to $[0,1],[1,0] \in \IP^1$.

Let $E_i$ be the general fiber of $\pi_i$, for $i = 1,2,3$. Then $\langle E_i, E_j \rangle = 2$ for any $i,j \in \{1,2,3\}, i \neq j$. 
The image of the morphism
$\pi_1 \times \pi_2 \times \pi_3 \colon X \rightarrow \IP^1 \times \IP^1 \times \IP^1$
is then defined by the tridegree $(2,2,2)$ polynomial
\[
 (u_{0}^2 - u_{1}^2)(v_{0}^2 - v_{1}^2) w_{0}^2 = f(u_{0},u_{1},v_{0},v_{1}) w_{1}^2.
\]
Consider the Segre embedding $\Sigma\colon \IP^1 \times \IP^1 \times \IP^1 \hookrightarrow \IP^7$, defined by
\begin{multline*}
 ([u_0,u_1],[v_0,v_1],[w_0,w_1]) \mapsto [x_0,x_1,x_2,x_3,x_4,x_5,x_6,x_7] = \\
= 
[u_0v_0w_0, u_0v_1w_1, u_1v_0w_1, u_1v_1w_0, u_0v_0w_1, u_0v_1w_0, u_1v_0w_0, u_1v_1w_1].
\end{multline*}
The involution on $\IP^7$ given by $[x_0,\ldots,x_7] \mapsto [x_0,\ldots,x_3,-x_4,\ldots,-x_7]$
induces the Enriques involution $\tilde\enr$ on $X$. Hence, we have the following commuting diagram
\[
 \begin{tikzcd}
 X \arrow[rr, "\pi_1 \times \pi_2 \times \pi_3"] \arrow[d] & & \IP^1 \times \IP^1 \times \IP^1 \arrow[r,hook,"\Sigma"] & \IP^7 \arrow[d, "\pr_{0123}"] \\
 Y \arrow[rrr] & & & \IP^3
 \end{tikzcd}
\]
where $\pr_{0123}$ is the projection $[x_0,x_1,x_2,x_3,x_4,x_5,x_6,x_7] \mapsto [x_0,x_1,x_2,x_3]$.
Note that the half-pencils are mapped onto the coordinate tetrahedron in $\IP^3$, so the image of $Y$ in $\IP^3$ is defined by an Enriques sextic.

\subsubsection{Nikulin-Kondo type \NKtype{1}} 
For $\alpha\in \IC \setminus \{ 1,\tfrac12, \tfrac32\}$, let $C$ be the curve defined by
\[
 C\colon (2\, u_{0}^2 - u_{1}^2)(v_{0}^2-v_{1}^2) = (2\,\alpha v_{0}^2 + (1 - 2\,\alpha)v_{1}^2)(u_{0}^2 - u_{1}^2).
\]
Put $B = C + \sum_{i = 1}^{4} L_i$. Then, the minimal resolution of the double covering 
of $\IP^1 \times \IP^1$ ramified over $B$ is a K3 surface $X$ endowed with 
an Enriques involution $\tilde\enr$ such that the quotient $X/\langle \tilde\enr \rangle$ has Nikulin-Kondo type \NKtype{1}.

Consider the curves
\begin{align*}
Q_{1}&\colon u_{1}v_{0} + u_{0} v_{1} = 0; \quad Q_2\colon u_{1}v_{0} - u_{0}v_{1} = 0; \\
Z&\colon \left(u_{0} + 3\, u_{1} \right) v_{0}^2 + \left( 3\, u_{0} + u_{1} \right) v_{1}^2 = 0.
\end{align*}
The curve $Z$ intersects $Q_{1}$ and $Q_2$ in one point with multiplicity $3$, and intersects~$C$ with even multiplicities if and only if 
\[
\alpha = \frac{15}{16} \quad \text{or} \quad \alpha = \frac{17}{16}.
\]
(The two cases differ only by a relabeling of the variables.)

In these cases, consider the sublattice $S' \subset S_X$ generated by the classes of the strict transforms of $C,L_{1},\ldots,L_4,Q_{1},Q_2,Z$ and of the exceptional divisors. 
Then, $\rank S' = 20$ and $\det S' = 7$, hence the same holds for $S_X$. This implies that $X$ is isomorphic to $X_{7}$, so the quotient $X/\langle \tilde\enr \rangle$ is isomorphic to $Y_\NKtype{1}$.

An Enriques sextic model for $Y_\NKtype{1}$ is given by
\begin{multline*}
(2 \, \alpha - 2) x_{0}^{2} x_{1}^{2} x_{2}^{2} 
+ x_{0}^{2} x_{1}^{2} x_{3}^{2} + x_{0}^{2} x_{2}^{2} x_{3}^{2} 
+ (2 \, \alpha -2) x_{1}^{2} x_{2}^{2} x_{3}^{2} = 
\\ = x_0 x_1 x_2 x_3 \left( x_{0}^{2} + ( 2 \, \alpha -3) x_{1}^{2} + ( 2 \, \alpha -1) x_{2}^{2} 
+ x_{3}^{2}\right).
\end{multline*}

\subsubsection{Nikulin-Kondo type \NKtype{2}} 
For $\alpha\in \IC \setminus \{0,-1\}$, let $C$ be the curve defined by
\[
 C\colon (v_{0}^2 - v_{1}^2)u_{0}^2 - (v_{0}^2 + \alpha v_{1}^2)u_{1}^2 = 0.
\]
Put $B = C + \sum_{i = 1}^{4} L_i$. Then, the minimal resolution of the double covering of $\IP^1 \times \IP^1$ ramified over $B$ is a K3 surface $X$ endowed with an Enriques involution $\tilde\enr$ such that the quotient $X/\langle \tilde\enr \rangle$ has Nikulin-Kondo type \NKtype{2}.

Consider the curves
\begin{align*}
F_{1}&\colon u_{1} = 0; \quad F_2\colon v_{1} = 0; \\
Z&\colon {\left(u_{0} - u_{1}\right)} v_{0} + {\left(u_{0} + 3 \, u_{1}\right)} v_{1} = 0
\end{align*}
The curve $Z$ intersects $C$ in a third point of multiplicity 2 exactly when 
\[
 \alpha = 63.
\]

In this case, consider the sublattice $S' \subset S_X$ generated by the classes of the strict transforms of $C,L_{1},\ldots,L_4,F_{1},F_2,Z$ and of the exceptional divisors. 
Then, $\rank S' = 20$ and $\det S' = 7$, hence the same holds for $S_X$. This implies that $X$ is isomorphic to $X_{7}$, so the quotient $X/\langle \tilde\enr \rangle$ is isomorphic to $Y_\NKtype{2}$.

An Enriques sextic model for $Y_\NKtype{2}$ is given by
\[
-x_{0}^{2} x_{1}^{2} x_{2}^{2} + x_{0}^{2} x_{1}^{2} x_{3}^{2} 
 + x_{0}^{2} x_{2}^{2} x_{3}^{2} +\alpha x_{1}^{2} x_{2}^{2} x_{3}^{2} 
= x_{0}x_1x_2x_3 \left( x_{0}^{2} - x_{1}^{2} - x_{2}^{2} + x_{3}^{2} \right).
\]

\subsection*{Acknowledgements} 
Both authors warmly thank Hisanori Ohashi and the other organizers of the 
3rd edition of the Japanese-European Symposium on Symplectic Varieties and Moduli Spaces at Tokyo University of Science in August 2018, where their collaboration started.

The first author would like to thank Igor Dolgachev, Shigeyuki Kondo and Shigeru Mukai
for  discussions.

The second author would like to thank Fabio Bernasconi, Chiara Camere, 
Alberto Cattaneo, Alex Degtyarev, Dino Festi, Grzegorz Kapustka, Roberto Laface and 
Matthias Schütt for their support and interest in this work. A special acknowledgement goes to Simon Brandhorst for his help with {\tt sage}.

\bibliographystyle{plain}
\bibliography{myrefsthemost}

\begin{thebibliography}{10}

\bibitem{AllcockDolgachev2018}
Daniel Allcock and Igor~V. Dolgachev.
\newblock The tetrahedron and automorphisms of {E}nriques and {C}oble surfaces
  of {H}essian type, 2018.
\newblock ar{X}iv:1809.07819.

\bibitem{bhpv}
Wolf Barth, Klaus Hulek, Chris Peters, and Antonius~Van de~Ven.
\newblock {\em Compact complex surfaces}, volume~4 of {\em Ergeb. Math.
  Grenzgeb.}
\newblock Springer-Verlag, Berlin Heidelberg New York, 2004.

\bibitem{women2015}
Marie~Jos\'{e} Bertin, Alice Garbagnati, Ruthi Hortsch, Odile Lecacheux, Makiko
  Mase, Cec\'{i}lia Salgado, and Ursula Whitcher.
\newblock Classifications of elliptic fibrations of a singular {K}3 surface.
\newblock In {\em Women in numbers {E}urope}, volume~2 of {\em Assoc. Women
  Math. Ser.}, pages 17--49. Springer, Cham, 2015.

\bibitem{Bor1}
Richard Borcherds.
\newblock Automorphism groups of {L}orentzian lattices.
\newblock {\em J. Algebra}, 111(1):133--153, 1987.

\bibitem{Bor2}
Richard Borcherds.
\newblock Coxeter groups, {L}orentzian lattices, and {$K3$} surfaces.
\newblock {\em Internat. Math. Res. Notices}, 1998(19):1011--1031, 1998.

\bibitem{Brandhorst}
Simon Brandhorst.
\newblock Private communication, 2018.

\bibitem{BS}
Simon Brandhorst and Ichiro Shimada.
\newblock {B}orcherds method for {E}nriques surfaces, 2018.
\newblock In preparation.

\bibitem{conway-sloane}
John~H. Conway and Neil J.~A. Sloane.
\newblock {\em Sphere Packings, Lattices and Groups}, volume 290 of {\em
  Grundlehren Math. Wiss.}
\newblock Springer-Verlag, Berlin Heidelberg New York, 1999.

\bibitem{Degtyarev}
Alex Degtyarev.
\newblock Private communication, 2018.

\bibitem{dolgachev-brief-introduction}
Igor~V. Dolgachev.
\newblock A brief introduction to {Enriques} surfaces.
\newblock In {\em Development of Moduli Theory -- Kyoto 2013, Adv. Study in
  Pure Math. Math. Soc.}, volume~69, pages 1--32. Math. Soc. Japan, 2016.

\bibitem{GAP}
The~GAP Group.
\newblock {G}{A}{P} - {G}roups, {A}lgorithms, and {P}rogramming.
\newblock Version 4.8.6; 2016 (\url{http://www.gap-system.org}).

\bibitem{harrache-lecacheux}
Titem Harrache and Odile Lecacheux.
\newblock Études des fibrations elliptiques d'une surface {K3}.
\newblock {\em J. Théor. Nombres Bordeaux}, 23:183--207, 2011.

\bibitem{horikawaI}
Eiji Horikawa.
\newblock On the periods of {Enriques} surfaces, {I}.
\newblock {\em Math. Ann.}, 234:73--108, 1978.

\bibitem{horikawaII}
Eiji Horikawa.
\newblock On the periods of {Enriques} surfaces, {II}.
\newblock {\em Math. Ann.}, 235:217--246, 1978.

\bibitem{HS2012}
Klaus Hulek and Matthias Sch\"utt.
\newblock Arithmetic of singular {E}nriques surfaces.
\newblock {\em Algebra Number Theory}, 6(2):195--230, 2012.

\bibitem{Keum1990}
JongHae Keum.
\newblock Every algebraic {K}ummer surface is the {$K3$}-cover of an {E}nriques
  surface.
\newblock {\em Nagoya Math. J.}, 118:99--110, 1990.

\bibitem{KeumKondo2001}
JongHae Keum and Shigeyuki Kondo.
\newblock The automorphism groups of {K}ummer surfaces associated with the
  product of two elliptic curves.
\newblock {\em Trans. Amer. Math. Soc.}, 353(4):1469--1487, 2001.

\bibitem{Kneser57}
Martin Kneser.
\newblock Klassenzahlen definiter quadratischer formen.
\newblock {\em Arch. Math.}, 8:241--250, 1957.

\bibitem{Kondo1986}
Shigeyuki Kondo.
\newblock Enriques surfaces with finite automorphism groups.
\newblock {\em Japan. J. Math. (N.S.)}, 12(2):191--282, 1986.

\bibitem{KondoKmJacC}
Shigeyuki Kondo.
\newblock The automorphism group of a generic {J}acobian {K}ummer surface.
\newblock {\em J. Algebraic Geom.}, 7(3):589--609, 1998.

\bibitem{lecacheux}
Odile Lecacheux.
\newblock Weierstrass equations for all elliptic fibrations on the modular
  {$K3$} surface associated to {$\Gamma_1(7)$}.
\newblock {\em Rocky Mountain J. Math.}, 45(5):1481--1509, 2015.

\bibitem{Lee2012}
Kwangwoo Lee.
\newblock Which {K}3 surfaces with {P}icard number 19 cover an {E}nriques
  surface.
\newblock {\em Bull. Korean Math. Soc.}, 49(1):213--222, 2012.

\bibitem{Mukai}
Shigeru Mukai.
\newblock Private communication, 2018.

\bibitem{MukaiOhashi2015}
Shigeru Mukai and Hisanori Ohashi.
\newblock The automorphism groups of {E}nriques surfaces covered by symmetric
  quartic surfaces.
\newblock In {\em Recent advances in algebraic geometry}, volume 417 of {\em
  London Math. Soc. Lecture Note Ser.}, pages 307--320. Cambridge Univ. Press,
  Cambridge, 2015.

\bibitem{Nikulin1979}
Viacheslav~V. Nikulin.
\newblock Integer symmetric bilinear forms and some of their geometric
  applications.
\newblock {\em Izv. Akad. Nauk SSSR Ser. Mat.}, 43(1):111--177, 238, 1979.
\newblock English translation: Math USSR-Izv. 14 (1979), no. 1, 103--167
  (1980).

\bibitem{Nikulin1984}
Viacheslav~V. Nikulin.
\newblock Description of automorphism groups of {E}nriques surfaces.
\newblock {\em Dokl. Akad. Nauk SSSR}, 277(6):1324--1327, 1984.
\newblock Soviet Math. Dokl. 30 (1984), No.1 282--285.

\bibitem{ohashi}
Hisanori Ohashi.
\newblock On the number of {Enriques} quotients of a {K3} surface.
\newblock {\em Publ. RIMS, Kyoto Univ.}, 43:181--200, 2007.

\bibitem{ohashi2}
Hisanori Ohashi.
\newblock Enriques surfaces covered by {J}acobian {K}ummer surfaces.
\newblock {\em Nagoya Math. J.}, 195:165--186, 2009.

\bibitem{Schutt2007}
Matthias Sch\"{u}tt.
\newblock Fields of definition of singular {$K3$} surfaces.
\newblock {\em Commun. Number Theory Phys.}, 1(2):307--321, 2007.

\bibitem{Sertoz2005}
Ali~Sinan Sert\"{o}z.
\newblock Which singular {$K3$} surfaces cover an {E}nriques surface.
\newblock {\em Proc. Amer. Math. Soc.}, 133(1):43--50, 2005.

\bibitem{ShimadaReduction}
Ichiro Shimada.
\newblock Transcendental lattices and supersingular reduction lattices of a
  singular {$K3$} surface.
\newblock {\em Trans. Amer. Math. Soc.}, 361(2):909--949, 2009.

\bibitem{ShimadaAlgo}
Ichiro Shimada.
\newblock An algorithm to compute automorphism groups of {$K3$} surfaces and an
  application to singular {$K3$} surfaces.
\newblock {\em Int. Math. Res. Not. IMRN}, 2015(22):11961--12014, 2015.

\bibitem{ShimadaSch}
Ichiro Shimada.
\newblock The automorphism groups of certain singular {$K3$} surfaces and an
  {E}nriques surface.
\newblock In {\em K3 surfaces and their moduli}, volume 315 of {\em Progr.
  Math.}, pages 297--343. Birkh\"{a}user/Springer, [Cham], 2016.

\bibitem{ShimadaHoles}
Ichiro Shimada.
\newblock Holes of the {L}eech lattice and the projective models of {$K3$}
  surfaces.
\newblock {\em Math. Proc. Cambridge Philos. Soc.}, 163(1):125--143, 2017.

\bibitem{thecompdata}
Ichiro Shimada.
\newblock Enriques involutions on singular {K3} surfaces of small
  discriminants: computational data, 2018.
\newblock \url{http://www.math.sci.hiroshima-u.ac.jp/~shimada/K3.html}.

\bibitem{Shimada2017}
Ichiro Shimada.
\newblock On an {E}nriques surface associated with a quartic {H}essian surface,
  2018.
\newblock ar{X}iv:1701.00580, to appear in \emph{Canad. J. Math.}

\bibitem{ShiodaInose}
Tetsuji Shioda and Hiroshi Inose.
\newblock On singular {$K3$} surfaces.
\newblock In {\em Complex analysis and algebraic geometry}, pages 119--136.
  Iwanami Shoten, Tokyo, 1977.

\bibitem{sagemath}
{The Sage Developers}.
\newblock {S}agemath, the {S}age {M}athematics {S}oftware {S}ystem ({V}ersion
  8.5), 2018.
\newblock \url{https://www.sagemath.org}.

\bibitem{Ujikawa2013}
Masashi Ujikawa.
\newblock The automorphism group of the singular {$K3$} surface of discriminant
  7.
\newblock {\em Comment. Math. Univ. St. Pauli}, 62(1):11--29, 2013.

\bibitem{Utsumi2016}
Kazuki Utsumi.
\newblock Jacobian fibrations on the singular {$K3$} surface of discriminant 3.
\newblock {\em J. Math. Soc. Japan}, 68(3):1133--1146, 2016.

\bibitem{Vinberg1975}
\`Ernest~B. Vinberg.
\newblock Some arithmetical discrete groups in {L}oba\v cevski\u\i \ spaces.
\newblock In {\em Discrete subgroups of {L}ie groups and applications to moduli
  ({I}nternat. {C}olloq., {B}ombay, 1973)}, pages 323--348. Oxford Univ. Press,
  Bombay, 1975.

\bibitem{Vinberg1983}
\`Ernest~B. Vinberg.
\newblock The two most algebraic {$K3$} surfaces.
\newblock {\em Math. Ann.}, 265(1):1--21, 1983.

\end{thebibliography}

\end{document}